\documentclass[a4paper]{amsart}
\usepackage{amsmath,amssymb,amsthm}

\usepackage{fancyhdr}
\usepackage{epsfig,amsfonts,latexsym}
\usepackage{color}

 \usepackage{url} 
 \usepackage{framed}
\usepackage{amsmath, color}
\usepackage{amssymb}
\usepackage{amssymb}
\usepackage{latexsym}
 \usepackage{graphicx, epstopdf, verbatim}
\usepackage{psfrag}

\newcommand{\R}{\mathbb R}
\newcommand{\Z}{\mathbb Z}

\newcommand{\HOX}[1]{}

\newcommand{\Afun}[1]{A(#1)}

\def\B{\mathcal B}

\newcommand{\figcommented}[1]{}

\def \beq {\begin {eqnarray}}
\def \eeq {\end {eqnarray}}
\def \ba {\begin {eqnarray*}}
\def \ea {\end  {eqnarray*}}

\renewcommand{\hat}{\widehat}
\renewcommand{\tilde}{\widetilde}

\newcommand{\ep}{\varepsilon}
\newcommand{\de}{\delta}

\newcommand{\co}{\colon}

\newcommand{\be}{\begin{equation}}
\newcommand{\ee}{\end{equation}}

\DeclareMathOperator{\dist}{dist}

\DeclareMathOperator{\diam}{diam}
\DeclareMathOperator{\inj}{inj}

\numberwithin{equation}{section}

\newtheorem{theorem}{Theorem}

\newtheorem{lemma}{Lemma}[section]

\newtheorem{proposition}[lemma]{Proposition}

\theoremstyle{definition}
\newtheorem{definition}[lemma]{Definition}

\theoremstyle{remark}
\newtheorem{remark}[lemma]{Remark}

\newcommand{\Sec}{\operatorname{Sec}}

\def\Prob{{\mathbb P}}
\def\Expec{{\mathbb E}}

\def\atext{}
\def\btext{} 
\def\btext{}
\def\trtext{}
\def\hhtext{}

\def\htext{}

\makeatother

 \def\Prob{\mathbb P}
 \def\Expec{\mathbb E}

\newcommand{\ignore}[1]{}

\newcommand{\barr}{\begin{array}}
\newcommand{\earr}{\end{array}}

\def\bfo{\begin{eqnarray*}}
\def\efo{\end{eqnarray*}}
\def\ba{\begin{eqnarray*}}
\def\ea{\end{eqnarray*}}
\def\beq{\begin{eqnarray}}
\def\eeq{\end{eqnarray}}

\def\hat{\widehat}
\def\tilde{\widetilde}

\renewcommand{\H}{{\mathbb H}} 
 
\newcommand{\I}{{\mathcal I}}

\newcommand{\proofbox}{\hfill $\square$\medskip}

\def\supp{\hbox{supp}\,}

\def\diam{\hbox{diam}}

\def\dist{\hbox{dist}}

\def\e{\varepsilon}
\def\p{\partial}

\def\tilde{\widetilde}

\def\picture #1 by #2 (#3){
   \vsquare to #2{
     \hrule width #1 height 0pt depth 0pt
     \hfill
     \special{picture #3}}}

\def\scaledpicture #1 by #2 (#3 scaled #4){{
   \dimen0=#1 \dimen1=#2
   \divide\dimen0 by 1000 \multiply\dimen0 by #4
   \divide\dimen1 by 1000 \multiply\dimen1 by #4
   \picture \dimen0 by \dimen1 (#3 scaled #4)}}

\def \B {{\mathcal B}}

\def \Z {{\mathbb {Z}}}

\def \R {{\mathbb {R}}}

\def \be {e}

\def \H2s {H^{s+1}_0(\partial M\times [0,T/2])}

\def \supp {\hbox{supp }}
\def \diam {\,\hbox{diam}\,}

\def \e {\varepsilon}

\def \pa0 {\partial _0}

\def \p {\partial}

\def\e{\varepsilon}

\def\tilde{\widetilde}

\def \mbeq {\begin {eqnarray}}
\def \meeq {\end {eqnarray}}

\def \B {{\mathcal B}}

\def \be {{ {e}}}

\def \vol {\hbox{vol}}

\title
[Reconstruction of a manifold]
{Reconstruction of a Riemannian manifold from noisy   intrinsic distances}

\author[Fefferman,  Ivanov, 
Lassas, Narayanan]{Charles Fefferman, Sergei Ivanov, \\ 
Matti Lassas, Hariharan Narayanan\\
$\ $\\ $\hbox{\rm To the memory of Yaroslav Kurylev.}\vspace{-5mm}$}

\address{\hspace{2cm} \linebreak
Charles Fefferman, Princeton University, Mathematics Department, USA.\hspace{7cm} 
\linebreak
 Sergei Ivanov,
St.~Petersburg Department of Steklov Institute of Mathematics,  Russia.
\hspace{6cm} 
\linebreak
Matti Lassas,
University of Helsinki,
 Finland.
\hspace{8cm} 
 \linebreak
Hariharan Narayanan, 
Tata Institute for Fundamental Research, 
 India. }

\date{} 

\begin{document}
\maketitle



\begin{abstract}
We consider reconstruction of a manifold, or,  invariant manifold learning,  where  a smooth Riemannian manifold $M$ is determined
from intrinsic  distances (that is, geodesic distances) of points in a discrete subset of $M$. In the studied problem the Riemannian manifold   $(M,g)$ is considered as an abstract metric space with intrinsic distances, not as an embedded submanifold of
an ambient Euclidean space. 
Let $\{X_1,X_2,\dots,X_N\}$ be
a set of $N$ sample points sampled randomly from
 an unknown Riemannian $M$ manifold. 
We assume that we are given  the numbers $D_{jk}=d_M(X_j,X_k)+\eta_{jk}$, where $j,k\in \{1,2,\dots,N\}$. Here, $d_M(X_j,X_k)$ are geodesic distances, $\eta_{jk}$ are independent, identically distributed random variables  such that $\mathbb E e^{|\eta_{jk}|}$ is finite. We show that
 when $N$ is large enough, it is possible to construct
an approximation of the Riemannian manifold $(M,g)$ with a large probability. This problem is a generalization of 
the geometric Whitney problem with  random measurement errors.
We consider also the case when the information on noisy distance $D_{jk}$ of points $X_j$ and $X_k$ is missing with some probability. In particular, we consider the case when we have no information on points that are far away.\end{abstract}

\noindent
{\bf Keywords:}   Inverse problems, Manifold learning, Geometric Whitney problem.

\section{Introduction}

Let  $M$ be a manifold and  $g$  an intrinsic Riemannian metric on it. Assume that one is given distances, $d_M(X_j,X_k)$, with random measurement errors, between points in a randomly sampled set $\{X_1,X_2,\dots,X_N\}$ of points of $M$. In this paper we ask,
how one can construct  a  Riemannian manifold $(M^*,g^*)$ from these data so that with a large probability, the distance (in Lipschitz-sense) of the constructed manifold 
$(M^*,g^*)$ to the original Riemannian manifold $(M,g)$ can be estimated.
The need of constructing the non-Euclidean, intrinsic metric is encountered in many applications, e.g.,\ in medical and seismic imaging, discussed in Section \ref{subsec: Applications of invariant}.

In the traditional manifold learning, for instance by using the ISOMAP algorithm introduced in the seminal paper  \cite{TSL}, one often aims to map points $X_j$ to points $Y_j=F(X_j)$ in an Euclidean space $\R^m$, where $m\geq n$ is as small as possible,
so that the Euclidean distances $\|Y_j-Y_k\|_{\R^m}$ are close to the intrinsic  distances $d_M(X_j,X_k)$ and find a submanifold $\tilde M\subset \R^m$ that is close
to the points $Y_j$. However, even in the ideal case when one is given an infinite set of points $X_j$ that form a dense subset of a smooth manifold $M$ and one has no measurement errors, 
finding a map $F:M\to \R^m$ for which the embedded manifold $F(M)=\tilde M\subset \R^m$ is isometric to $(M,g)$  is numerically a very difficult 
task as it means finding a map which existence is proved by the  Nash embedding theorem,
see \cite{Nash1,Nash2} and \cite{Verma1} on numerical techniques based on  Nash embedding theorem. 
%
%
One can  overcome this difficulty by formulating the problem in a coordinate invariant way:
Given the geodesic distances of points sampled from a Riemannian manifold $(M,g)$, construct a manifold $M^*$ with an intrinsic metric tensor $g^*$ so that Lipschitz distance of $(M^*,g^*)$  to the original manifold $(M,g)$ is small. This problem was studied in  \cite{Meila,Meila2}
using diffusion maps \cite{diffusion,CoifmanLafon}. In this paper we 
consider this problem when distances  are given with random errors and
use metric geometry to construct $(M^*,g^*)$ so that the distance of
$(M^*,g^*)$ and $(M,g)$  
 can be estimated with a large probability. We emphasise that we consider 
$M^*$ as an abstract manifold, that is not isometricly embedded to an Euclidean space, but where the metric is given 
by a metric tensor $g^*$ that is constructed from the above data.

\subsection{The main result}

Let $n\ge 2$ be an integer, $\Lambda>0,$  $D>0$, and $i_0>0$. Let
$(M,g)$ be a compact Riemannian manifold of dimension $n$
such that
\beq \label{26.1} 
& & i)\  \| \hbox{Sec}_M\|_{L^\infty(M)}\leq \Lambda^2, 
\quad  ii)\  \diam (M) \leq D, \quad iii) \ \hbox{inj}\, (M) \geq i_0,
\eeq
where $\hbox{Sec}_M$  is the sectional curvature  of $(M,g)$, $\diam(M)$  is the diameter of $(M,g)$ 
and $\hbox{inj}\, (M)$  is the injectivity radius of $(M,g)$, that is, the minimal radius of Riemannian normal coordinates. 
Let $d_M(x,y)$
denote the intrinsic (or geodesic) distance of the points $x,y\in M$ determined by the metric tensor $g$
corresponding to the line element $ds^2=g_{jk}(x)dx^jdx^k$. Here and below, we use Einstein's summation convention and sum over indexes appearing as super and sub-indexes.

Let $(\Omega,\Sigma,\Prob)$ be a complete probability space,
$\B$ be the $\sigma$-algebra of Borel  sets on $M$, and
 $\mu:\B\to [0,1]$  be a probability measure on $M$.
Let $dV_g$  be Riemannian volume on $(M,g)$.
Assume that the Radon-Nikodym derivative of $\mu$  satisfies
\beq\label{Radon Nicodym}
0<\rho_{min}\leq \frac {d\mu}{d V_g}\leq \rho_{max},\quad\hbox{where $\rho_{min},\rho_{max}\in \R_+$.}
\eeq

\begin{definition}\label{def: X and D} Let 
$X_j$, $j=1,2,\dots,N$ be independent, identically distributed (i.i.d.) random variables having distribution $\mu$.
{\btext Let $\sigma>0,$ $\beta>1$, and $\eta_{jk}$ be random variables satisfying
\beq\label{noise 1}
\Expec \eta_{jk}=0,\quad 
\Expec ({\eta_{jk}^2})={\sigma^2},\quad \Expec e^{|\eta_{jk}|}= {\beta}.
\eeq
We assume that  all random variables $\eta_{jk}$ and $X_j$ are independent.} Let
\beq\label{data 1}
D_{jk}=d_M(X_j,X_k)+\eta_{jk}.
\eeq
be the geodesic distances of points $X_j$ and $X_k$ measured with errors $\eta_{jk}$.
\end{definition}

{\btext Note that the above assumptions are satisfied when  $\eta_{jk}\sim N(0,\sigma^2)$  are i.i.d.\ Gaussian random variables and $\beta\leq 2e^{\sigma^2}$.}
We are mostly interested in a case when $\sigma$ is fixed and $N$ is large.

{\hhtext
\begin{definition} \label{def 1}The partial data
 is given by
\ba
\overline D_{jk}=D_{jk}^{\hbox{\tiny (partial\ data)}}=\left\{\begin{array}{cl}  D_{jk} &\hbox{if $Y_{jk}=1$,}\\
 \hbox{`missing'} &\hbox{if $Y_{jk}=0$,}\end{array}\right.
\ea
where $Y_{jk}$  are random variables taking values in $\{0,1\}$ and $j,k\in \{1,2,\dots,N\}$.
We assume that $Y_{jk}$ are independent of random variables  $X_{j'}$ for all $j' \in\{1,2,\dots,N\}\setminus \{j,k\}$
and of  $\eta_{j''k''}$ for all $j''$  and $k''$. Above, the value `missing', can be replaced by a large real value, e.g. by  $2D$.

Assume that the
conditional probability of the event $\{Y_{jk}=1\}$, when $X_j$ and $X_k$  are known, is
\beq\label{Phi cond 1}
\Prob(Y_{jk}=1\,|\, X_j, X_k)=\Phi(X_j,X_k).
\eeq
More precisely, when 
 $\B_{jk}\subset \Sigma$ is the 
$\sigma$-algebra generated by the random variables 
$X_{j}$ and $X_k$, above in formula  (\ref{Phi cond 1}) we use notation $\Prob(Y_{jk}=1\,|\, X_j, X_k)
=\Prob(Y_{jk}=1\,|\, \mathcal B_{jk})$.
Here, $ \Phi:M\times M\to [0,1]$ is a measurable function such that
there is  a  function
$ \Phi^1:[0,\infty)\to [0,1]$ so that
$s\mapsto \Phi^1(s)$ is non-increasing and
 \beq\label{Phi cond 2}
\Phi^1(0)=\phi_0,\quad \|\Phi^1\|_{C^{1}(\R)}\leq H,\quad c_1\Phi^1(d_M(x,y))\leq 
\Phi(x,y) \leq c_2\Phi^1(d_M(x,y)),\quad x,y\in M,\hspace{-1cm}
\eeq
where $0<c_1<1<c_2$ and  $\phi_0\in \R_+$.  

\end{definition}

Below, for $t\in \R$ we denote by $\lfloor t \rfloor$  the largest integer $m$ such that
$m\leq t$.

We will show that probabilistic considerations involving the above data, combined with the deterministic results in \cite{FIKLN} (where we considered small deterministic errors) and Appendix A, yield that one can construct a smooth manifold $(M^*,g^*)$
that approximates the original manifold $(M,g)$. The proofs of Theorem \ref{thm 2:manifold}  below and the 
 results in \cite{FIKLN} give a procedure, which the output is a submanifold
 $M^*\subset \R^d$ (where $d$    depends only on $n,D,\Lambda,$ and $i_0$) and a metric tensor $g^*$ on $M^*$.

 \begin{theorem}
\label{thm 2:manifold}
%

{\atext Let $n\ge 2$, $D,\Lambda,i_0,\rho_{min},\rho_{max},\sigma,\beta,c_1,c_2,H, \phi_0>0$ be given.
Then there are $\delta_0>0$, and 
 $C_0>0$, depending on $n,D,\Lambda,i_0,\rho_{min},\rho_{max},\sigma,\beta,c_1,c_2, H, \phi_0$, and  there is ${C_1}>0$,
 depending on $n$,
  such that 
the following holds for  $\theta\in (0,\frac 12)$. }

Let   $M$  be a compact $n$-dimensional manifold satisfying (\ref{26.1}),   $0<{{\delta}}<\delta_0$, and
$$
N=\bigg\lfloor
C_0  \delta^{{ -3n}}
 \bigg(\log^2 (\frac{1}{\theta})+\log^8 (\frac 1 {{{\delta}}})\bigg)\bigg\rfloor 
$$
 and $\overline D_{jk}$, $j,k=1,2,\dots,N$ be as in Definitions \ref{def: X and D} and \ref{def 1}. 
Suppose that one is given 
samples of the random variables $\overline D_{jk}$  for $j,k=1,2,\dots,N$. Then
 withß a probability larger than $1-\theta$ one can 
 construct a  compact, smooth $n$-dimensional  Riemannian manifold $(M^*,g^*)$ that approximates the manifold $(M,g)$ in the following way:

\begin{enumerate}
\item 
There is a diffeomorphism $F:M^*\to M$ satisfying 
\beq
\label{Lip-condition}
\frac 1L\leq \frac{d_M(F(x),F(y))}{d_{M^*}(x,y)}\leq L,\quad \hbox{for all }x,y\in {M^*},
\eeq
where 
$L=1+C_1{{\delta}}$, that is, the Lipschitz distance 
 of the metric spaces $(M^*,g^*)$ and 
$(M,g)$ satisfies  $d_{Lip}((M^*,g^*),(M,g))\leq \log L$.

\item The sectional curvature $\hbox{Sec}_{M^*}$ of ${M^*}$ satisfies $|\hbox{Sec}_{M^*}|\le C_1\Lambda^2$.

\item The injectivity radius $\hbox{inj}({M^*})$ of ${M^*}$ satisfies
$$
 \hbox{inj}({M^*})\ge \min\{ {\atext(C_1\Lambda)^{-1}}, (1-C_1{{\delta}})\hbox{inj}(M)\} .
$$
\end{enumerate}
\end{theorem}

We note that the knowledge of the authors, the results are new also in the case when there is no missing data,
that is, $\Phi(x,y)=1$ for all $x,y\in M$. 
\medskip

\begin{remark} Theorem \ref{thm 2:manifold} concerns the regime where the noise level $\sigma$ is a fixed constant, the number of points
$N$ is large, and we are interested in the situation where we want the probability $\theta$ of a wrong final reconstruction to be very small. This is reflected by
the fact that the probability $\theta$  of obtaining a wrong reconstruction appears only in the logarithmic term $\log (\theta^{-1})$.
\end{remark}

\subsection{Idea of the proof and three nets of points on the manifold}

Let us assume that $N=N_0+N_1+N_2$, where $N_0,N_1,N_2\in\Z_+$. We are interested in the case when
$N_2>N_1>N_0$.
We consider the random set $S_0=\{X_1,\dots,X_{N_0}\}$ as a coarse net on $M$ and
compute approximate distances between the points in the net $S_0$ by using the auxiliary nets
$S_1=\{X_{N_0+1},\dots,X_{N_0+N_1}\}$ and $S_2=\{X_{N_0+N_1+1},\dots,X_{N_0+N_1+N_2}\}$. 
These random sets correspond to the index sets
\ba
I^{(0)}=\{1,2,\dots,N_0\},\quad I^{(1)}=\{N_0+1,\dots,N_0+N_1\},\quad I^{(2)}=\{N_0+N_1+1,\dots,N_0+N_1+N_2\}.
\ea
Recall that $X_j$ are independent random variables, taking values on $M$, with the distribution $\mu$.

Below, we say that a set $Y\subset M$  is $\delta-$dense in $M$ if for all $p\in M$  there is $y\in Y$
such such that $d_M(p,y)<\delta$. A  $\delta-$dense subset $S$ of $M$ is often called a $\delta-$net.


 
Let us give an overview of the ideas on the proof of Theorem \ref{thm 2:manifold}. 

First, we  use the ``densest net'' $S_2$  to compute approximately the numbers
\beq\label{k Phi function}
k_\Phi(y,z)=\int_M |d_M(y,x)-d_M(z,x)|^2\,\Phi(y,x)\Phi(x,z)\, d\mu(x),
\eeq	
for $y$  and $z$ in the ``medium dense net'' $S_1$, see Prop. \ref{prop: L2 norm of distance functions}.
This corresponds to taking average of function $|d_M(y,x)-d_M(z,x)|^2$ over all those
sample points $x\in S_2$ for which the distances $d_M(y,x)$  and $d_M(z,x)$  are not missing, see \eqref{Kjk2 B}.

Note that when the product $\Phi(y,x)\Phi(x,z)$  is small, there are only a small amount
of sample points  $x\in S_2$  for which the value of the function $|d_M(y,x)-d_M(z,x)|^2$ 
can be computed, and then the estimator for the function $k_\Phi(y,z)$ is not reliable.
Thus the reliability of the estimator for the function $k_\Phi(y,z)$ is measured
by
\beq
A_\Phi(y,z)&=&\int_M \Phi(y,x)\Phi(x,z)\, d\mu(x).
\eeq	
Indeed, when $A_\Phi(y,z)$ is larger than some threshold value $b>0$, the obtained estimator
for the function $k_\Phi(y,z)$ is reliable with a large probability. When we compute 
an estimator for the function $k_\Phi(y,z)$ using a sampling imitating the integral in \eqref{k Phi function},
we can compute also an estimator for $A_\Phi(y,z)$, see \eqref{Kjk2}.

\centerline{\includegraphics[width=3cm]{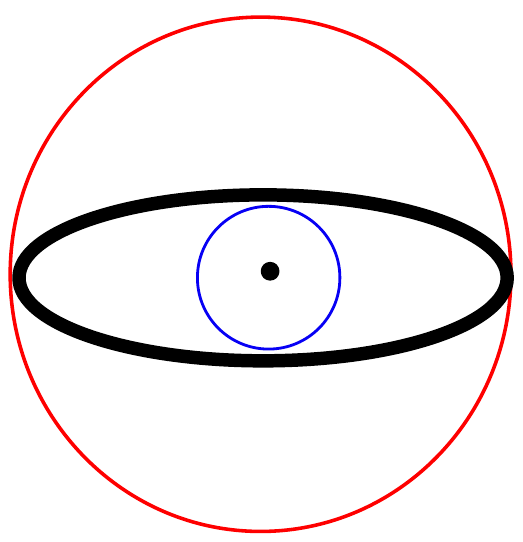}
\qquad
\includegraphics[width=5cm]{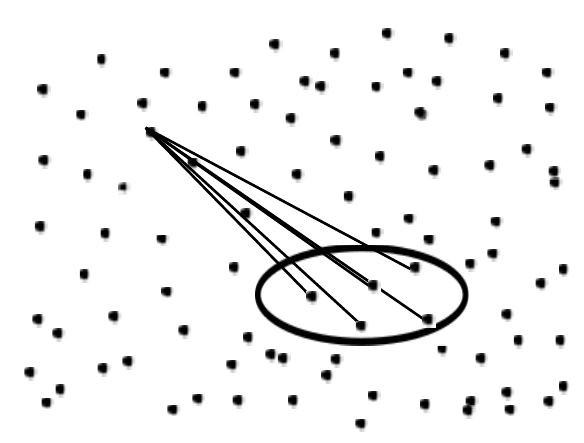}}

\noindent {{\bf Figure 1.} \it 
Left:  
The set 
$D_\Phi(y_1,\rho)$
 satisfies $ B_M(y,c_5^{-1}\rho)\subset D_\Phi(y_1,\rho)\subset B_M(y,\rho)$ and thus
 $D_\Phi(y_1,\rho)$ can be considered as an approximate $\rho$-neighbourhood of the the point $y$.
 Right: The approximate distance 
$d^{app}(y_1,y_2)$ in the formula (\ref{d1 app CC}) is the average of distances from $y_2$ to the points in the neighbourhood $D_\Phi(y_1,\rho)$.
Later, we approximate $d^{app}(y_1,y_2)$ by taking the average of distances of $y_2$ to the points in $S_1\cap D_\Phi(y_1,\rho)$, where $S_1$
is the medium dense net of sample points.
}
\medskip

Second, we are going to use the  set $S_1$  to compute the 
approximate distances $d^{app}(y_1,y_2)$  of the points $y_1$ and $y_2$ in
the ``coarse net'' $S_0$ using reliable distances  $k_\Phi(y,z)^{1/2}$. We do this by computing estimators
for the functions (see  Definition \ref{Q term 1})
\beq\label{d1 app BB missing}
& &Q(y_1,y_2)=\frac {V_\Phi(y_1,y_2)}{W_\Phi(y_1,y_2)},\\
\nonumber 
& &V_\Phi(y_1,y_2)= \int_{M} 
\beta_1(\frac {A_\Phi(y_1,z)}b)\, 
\psi_\rho(k_\Phi(y_1,z))
\Phi(z,y_2)d_M(z,y_2)d\mu(z),\\
\nonumber 
& &W_\Phi(y_1,y_2)= \int_{M} 
\beta_1(\frac {A_\Phi(y_1,z)}b)\psi_\rho(k_\Phi(y_1,z))\Phi(z,y_2)
d\mu(z),
\eeq
where $\beta_1\in C^\infty(\R)$ is a cut-off function such that
$\beta_1(t)=0$ for $t<1$ and $\beta_1(t)=1$ for $t>2$, 
and $\psi_\rho \in C^\infty(\R)$ is a cut-off function such that
$\psi_\rho(s)=1$ for $s<\rho^2$ and $\psi_\rho(s)=0$ for $s>2\rho^2$.
Here, $\beta_1({A_\Phi(y_1,z)}/b) \psi_\rho(k_\Phi(y,z))$ is the smoothened version of the indicator function of
the set 
\beq\label{DAP}
D_\Phi(y_1,\rho)=\{z\in M:\ k_\Phi(y,z)<\rho^2,\ A_\Phi(y_1,z)\geq b\}.
\eeq
The set is $D_\Phi(y_1,\rho)$ is a Lipschitz approximation the union 
of the  ball $B_M(y_1,\rho)$.

Then, roughly speaking, we compute an estimator for the function  $V_\Phi(y_1,y_2)$ 
computing averages of function
$\beta_1({A_\Phi(y_1,z)}/b)d_M(z,y_2)$
over all sample points $z$ the medium net  $S_1$ that are in the set
$D_\Phi(y_1,\rho)$ and for which data on the distance $d_M(z,y_2)$
is not missing.
At the same time, we compute an estimator for the function  $W_\Phi(y_1,y_2)$ by
computing averages of function
$\beta_1({A_\Phi(y_1,z)}/b)\psi_\rho(k_\Phi(y,z))$
over the same sample points.
The idea is that when $W_\Phi(y_1,y_2)$ is larger than some threshold $u>0$, the 
estimators computed from random data for the functions $V_\Phi(y_1,y_2)$, $W_\Phi(y_1,y_2)$, and $Q(y_1,y_2)$ 
are reliable with a large probability. 
Then, we define for $y_1,y_2\in S_0$
\beq\label{d1 app CC}
& &d^{app}(y_1,y_2)=\left\{\begin{array}{cl}{Q}(y_1,y_2), &\hbox{if $W_\Phi(y_1,y_2)>u
$},\\
 D, &\hbox{otherwise.}\\ \end{array}\right.
\eeq
Then there is $r_1$ such that
$d^{app}(y_1,y_2)$ approximate the true distance $d_M(y_1,y_2)$ with a small error $\e_1$ when $d_M(y_1,y_2)<r_1$,
and moreover, if $d_M(y_1,y_2)\ge r_1$, then $d^{app}(y_1,y_2)>r_1-\e_1$. 
In other words, with a large probability,
we can construct the distances  $d_M(y_1,y_2)$
with small errors for all points
$y_1$ and $y_2$ in $S_0$ that close to each other.

After the above constructions, we will use Proposition \ref{prop:improved density} in Appendix A, concerning a reconstruction of a Riemannian manifold
when we are given distances with small (deterministic) errors. This result is an improved version of the earlier results given in \cite[Corollary 1.10]{FIKLN}.

\subsection{Earlier results for submanifolds of $\R^n$ and graphs of functions.}\label{subsec: Earlier}



%
%
%

In dimensionality reduction and in the traditional manifold learning, the aim is to  transform data, consisting of points in a $d$-dimensional space that are   near an $n$-dimensional submanifold $M$, where $d>>n$ 
into a set of points in a low dimensional space $\R^m$ close to an $n$-dimensional
submanifold, where $d>m\ge n$. During transformation all of them try to preserve some geometric
properties, such as appropriately measured distances between points of the original data set, see
\cite{cheng,CB,TSL}. 
Perhaps the most basic of such methods is `Principal Component Analysis' (PCA), \cite{PCA1,kernel}
where 
one projects the data points onto the span of the $n$  eigenvectors corresponding
to the $n$ largest eigenvalues of the ($d\times d$) covariance matrix of the data points.


In the case of `Multi Dimensional Scaling' (MDS) \cite{mds}, the pairwise distances between
points are attempted to be preserved. One minimizes a certain `stress function' which captures the total error in pairwise distances between the data points
and between their lower  dimensional counterparts. For instance, given points $(x_j)_{j=1}^N$, $x_j\in \R^d$,
one tries to find $(y_j)_{j=1}^N$, $x_j\in \R^m$, that
is an (approximate)  minimizer of a stress function
 \beq\label{MDS minimization problem}
\min_{y_j\in \R^m} (\sum_{i,j=1}^N (\|y_i-y_j\|_{\R^m}-d_{ij})^2),
\eeq 
where $d_{ij}=\|x_i-x_j\|_{\R^d}$ are the Euclidean distances of points $x_i$ and $x_j$.

`ISOMAP' \cite{TSL} attempts to improve on MDS by trying to capture geodesic distances
between points while projecting. For each data point $x_i$ in the data set $\mathcal X=(x_j)_{j=1}^N$,  $x_j\in \R^d$, a `neighbourhood graph' 
is constructed using the  
$K$-neighbours of $x_i$, that is, the $K$ nearest points of $\mathcal X$  to $x_i$, the edges 
carrying the length between points. Now the shortest distance between points is computed in the resulting global graph containing all the neighbourhood graphs 
using a standard graph theoretic algorithm such as Dijkstra's.
Let $D ^G= [d^G_{ij}]$ be the $N\times N$ matrix of graph distances. 
Then MDS is used to find $(y_j)_{j=1}^N$, $x_j\in \R^m$ that (approximately) solve the minimization problem 
(\ref{MDS minimization problem}) with distances $d_{ij}$ replaced by $d^G_{ij}$.
If the data set $\mathcal X=(x_j)_{j=1}^N$ consists of $\delta$-dense set points of a submanifold $M\subset \R^d$ with small $\delta$,
then ISOMAP tries to find an approximation for isometric embedding, that is, a map $F:M\to \R^m$ for which 
 \beq\label{isometry minimization problem}
\|F(x)-F(y)\|_{\R^m}\approx d_M(x,y),\quad x,y\in M,
\eeq
where $d_M(x,y)$ is the intrinsic distance of the points $x$ and $y$ of the isometrically embedded manifold $M\subset \R^d$. 
Observe that one can not have equality $\|F(x)-F(y)\|_{\R^m}= d_M(x,y)$ in (\ref{isometry minimization problem})
unless all unit speed geodesics of $M$  are mapped to Euclidean lines (parametrized by the Euclidean length) in the map $F$, implying that $M$ has to be a flat manifold.
Thus, even if the data points are contained in a submanifold $M$ and the number $N$ of data points grows and if $d_{ij}$ are equal to the 
intrinsic distances $d_M(x_i,x_j)$, the ISOMAP algorithm does not produce a manifold which intrinsic distances are the same as those of the original manifold $M$ when the manifold  $M$ has non-zero curvature. The convexity and flatness conditions that guarantee that  the ISOMAP algorithm reconstructs the original manifold are studied in \cite{Bernstien,Donoho1,Donoho} for ISOMAP and in \cite{Zha} for the continuum ISOMAP. 

%

In the seminal papers \cite{diffusion,CoifmanLafon} on `Diffusion Maps', 
a complete graph is built on the data points sampled from a manifold $(M,g)$  and each
edge is assigned a weight $a(x,y)$ that is a Gaussian function of the distance of the points
$x$ and $y$. The normalized version of the kernel $a(x,y)$ defines a diffusion operator
on $M$ and this operator has eigenfunctions $\phi_j(x)$. These functions can be used to construct
a non-isometric embedding $x\mapsto (\phi_j(x))_{j=1}^m$ of the manifold $M$ into $\R^m$.
This construction  is continued in \cite{Meila} by computing an approximation the metric tensor $g$  by
using finite differences to find the Laplacian of  the products of the local coordinate functions. 
When there are no errors in the data, this construction is shown in \cite{Meila2} to converge to a correct limit as the number of  points sampled from the manifold $M$ tends to infinity. 
Other topological embedding methods for manifolds, based on heat kernels and eigenfunctions,
have been developed in \cite{BN,donoho,heat}. 
Moreover, locally linear construction methods are studied in \cite{mc,Riemannian,Riemannian2,llc,ZZ}.

The construction of a surface approximating a set of points in $\R^m$ is closely related to 
 the classical Whitney's problem. This problem is the construction
of a function $  F(x) \in C^m(\overline S),$ where $S \subset    \R^n$ is open, which is equal to a given function
$f(x)  $ on $  K$, where $K\subset S$.
This problem has been studied in different norms in \cite{F1,F2,F3,FK1} and the interpolation results on Whitney problem have been applied
for manifold of a submanifold of $\R^m$ in \cite{putative,FMN}.

%
%



\subsection{Applications of new results with intrinsic distances}\label{subsec: Applications of invariant}
{\it A. Submanifolds of Banach spaces.} In several imaging problems the images are considered as elements of Banach spaces that have no inner product structure.
For example, in many applications the 2-dimensional images (corresponding photographs) or  3-dimensional objects, are modeled by functions $u:D\to \R$, $D\subset \R^m$, $m=2,3$ that are
in
 the space of functions of bounded variation, $BV(D)$ or in Besov spaces $B^s_{11}(D)$, \cite{OsherRudin,LSS}. For example, in medical imaging,
 such functions are used to model  piecewise constant or piecewise smooth functions that correspond to the structure of the human body with internal organs with sharp jumps of density at the boundaries of the organs. 

{\it B. Physical models with  non-Euclidean metric.} 
Many inverse problems can be formulated as geometric problems
where the goal is to determine the underlying manifold structure.
For example, by probing a
medium with waves one can measure the travel times between points. This defines a non-Euclidean metric called the travel
time metric $g$. Recovering the wave speed function inside the medium is equivalent
to the determination of a Riemannian manifold from external or boundary measurements. The relation of the boundary measurements to the distances between the points in 
an $\delta$-net in the interior of the manifold is  considered e.g.\ in \cite{AKKLT,HoopOksanen,KKL,KOP}. 
Determining the wave speed of elastic waves inside a body  is a central problem in seismic
imaging of the Earth, \cite{SUV,U}.
An  example, in medical imaging, where the physical structures are represented using abstract Riemannian manifolds is ultrasound imaging where the acoustic properties of the inside of a body is imaged.
The typical ultrasound images correspond to image of the body represented in the non-Euclidean travel time coordinates,
of more precisely, in the Riemannian normal coordinates of the travel time metric \cite{Qu}. In these coordinates, the image rays, that is, the geodesics from the location of the source device are 
straight lines.
Theorem \ref{thm 2:manifold} can be applied when we  can measure travel times of waves between points. For instance, for the earthquakes 
the travel times of the surface waves can be directly measured and one can use these data to deduce
the properties of the Earth close to the surface. 
Another example is the Magnetic Resonance Imaging (MRI) based Elastography in medical imaging, where the elastic properties of the body of a patient are determined by observing the propagation of elastic waves sent into the body \cite{Hoskins}. 

Models involving local and missing data  are encountered in sensor technology, e.g. in the 
radio frequency identification (RFID) or in smart dust sensors, where a large number of low quality sensor send signals, either to receivers or to each others.
On the earlier results on missing data in manifold learning, see e.g. \cite{CL,GS}.
%
%
%
%
%
%

\subsection{Notations}

{\btext 
For a  Riemannian manifold $(M,g)$ satisfying (\ref{26.1}), let $\hbox{vol}_g$ be the Riemannian volume on $(M,g)$, and $p\in M$, the exponential map $\exp_p:T_pM\to M$
defines a smooth surjective map $\exp_p:\{\xi\in T_pM:\ \|\xi\|_g<D+1\}\to M$. {\atext By Bishop-Gromov inequality \cite[Chapter 9, Lem.\ 1.6]{Pe},
the function $r\mapsto \hbox{vol}_g(B(x,r))/v(n,-\Lambda^2,r)$  is non-increasing and bounded by 1, where $\hbox{vol}_g(B(x,r))$  
is the volume of the ball $B(x,r)\subset M$ and $v(n,-\Lambda^2,R)$ is the volume of the ball of radius $r$ in the hyperbolic space of dimension $n$ having 
curvature $-\Lambda^2$. Hence
\beq \label{26.1 B}
\vol_g(M)\leq V_0=v(n,-\Lambda^2,D)\leq {\omega_n} \bigg(\frac {\sinh(\Lambda D)}{\Lambda D}\bigg)^{n-1} D^n
\eeq
where $\omega_n$ is the volume of the unit ball in $\R^n$, see  \cite[Ch. 6, Cor. 2.4]{Pe}. Moreover,
\beq\label{estimate for number in ball 0}
\frac {\hbox{vol}_g(B(x,\rho))}{\hbox{vol}_g(M)} \geq  \frac {v(n,-\Lambda^2,\rho)}{v(n,-\Lambda^2,D)}
\geq 
\hat c_3\rho^n,\quad 
 {\mu(B(x,\rho))} \geq c_3\rho^n,\quad \hat c_3=\frac {\omega_n} {V_0},\quad c_3=\frac{\rho_{min}} {\rho_{max}}\hat c_3,\hspace{-2cm}
\eeq
where we use the fact that $v(n,-\Lambda^2,\rho)\ge \omega_n\rho^n$.
 Let
 \beq\label{r0 and r1 and phi1}
  r_0=\min(\frac 1{2H}\phi_0,i_0,\frac\pi{2\Lambda}),\quad r_1=\frac 12 r_0,\quad  \phi_1=\frac {c_1}2\phi_0.
\eeq
Then} Definition \ref{def 1} implies 
\beq\label{inequality phi 1} 
\Phi(x,y)\geq c_1\Phi^1(r_0)\geq \phi_1,\quad \hbox{for $d(x,y)\le r_0$}.
\eeq

\noindent{\bf Acknowledgements.}  The authors express their gratitude to  the Fields Institute,
where parts of this work have been done.
Ch.F.\  was partly supported {by} AFOSR, grant   DMS-1265524, and NSF, grant FA9550-12-1-0425.
S.I.\ was partly supported {by} RFBR, {grants 14-01-00062 and 17-01-00128-A}, 
M.L.\ was  supported by AF, grants  284715, 312110. 
 and H.N. was partly supported by NSF grant DMS-1620102 and a Ramanujan Fellowship.

%
\section{Reformulation of the main result with several parameters}

Our aim is to prove the following result that yields of Theorem \ref{thm 2:manifold} when it is combined with the results
in Appendix A
on manifold reconstruction with small deterministic errors.

\begin{theorem}\label{proposition Q final introduction}
Let $n\ge 2$  and $D,\Lambda,i_0,\rho_{min},\rho_{max},c_0,c_1,H,\sigma,\beta>0$ be given.
Then 
there are
 $C_2>1$, $\hat \e_1<1$, and $\hat \delta_1<1$ depending on $n,D,\Lambda,i_0,\rho_{min},\rho_{max},c_0,c_1,H,\sigma,\beta$,
 such that the following holds for  $\e_1\leq \hat \e_1,\delta_1\leq \hat \delta_1,$ and $\theta \in (0,1/2)$: 
 Let
\beq\label{N2 condition}
& &  N_0=\lfloor  C_2  \delta_1^{-n} ( {\log(\frac 1\theta)}+\log(\frac 1{\delta_1}))\rfloor,\quad N\ge N_0+\lfloor
C_2 \e_1^{-2n} 
 \bigg(\log^2 (\frac{1}{\theta})+\log^2(\frac 1 {\delta_1})+\log^8 (\frac 1 {\e_1})\bigg)\rfloor .
\eeq

Also, let $X_j$ and $\overline D_{jk}$, $j,k=1,2,\dots,N$ be as in Definitions \ref{def: X and D} and \ref{def 1}.
Suppose that we are given 
samples of the random variables $\overline D_{jk}$  for $j,k=1,2,\dots,N$. Let $r_1$ be given in (\ref{r0 and r1 and phi1}). Then
 with a probability larger than $1-\theta$ the set $\{X_j:\ j=1,2,\dots,N_0\}$
 is a $\delta_1$-net in $M$  and 
 one can
determine the approximate distances $d^{(a)}(X_j,X_{j'})$ so that the following holds:


For all $j,j'\in \{1,2,\dots,N_0\}$,
\beq\label{eq: comparison of approximate distances 1}
& &|d^{(a)}(X_j,X_{j'})-d_M(X_j,X_{j'})|\leq \e_1,\quad \hbox{if }d_M(X_j,X_{j'})<r_1,\\
\label{eq: comparison of approximate distances 2}
& &d^{(a)}(X_j,X_{j'})\geq r_1- \e_1,\quad \hbox{if }d_M(X_j,X_{j'})\ge r_1.
\eeq

In the case when the probability $\Phi(x,y)$, for that the information on the distance of $x$ and $y$ is not missing, is bounded from below by a positive constant, that is, $\Phi(x,y)\geq \phi_2,$
 for  all $(x,y)\in M\times M$, the inequality (\ref{eq: comparison of approximate distances 1}) holds for all $X_j,X_{j'}$,
$j,j'\in \{1,2,\dots,N_0\}$.
%
\end{theorem}

We note when $\Phi(x,y)$ is bounded from below by a positive constant, $\Phi(x,y)\ge c_1\phi_0,$
 we can choose
the function $\Phi^1$ to be equal to the constant $\phi_0$. Without loss of generality,
we can assume that in this case $\phi_2=c_1\phi_0$.

\subsection{Probability that the sample points form a dense net} 
First we estimate the probability that the set $S_0=\{X_1,\dots,X_{N_0}\}$ is a $\delta_1$-net, using standard methods based on  the collectors problem.

\begin{lemma}\label{Lemma dense net}
There is 
 $C_3
 \ge 10$ such if $\theta\in (0,\frac 12)$, $\delta_1\in (0,\frac 
D2)$, and
\beq\label{e: N0}
 N_0\geq  C_3  \delta_1^{-n} (\log(\delta_1^{-1})+ {\log(\theta^{-1})}),
\eeq
then the probability that the set $\{X_j:\ j=1,2,\dots,N_0\}$ is a $\delta_1$-net in $M$ is larger than $1-\frac 12\theta$.
\end{lemma}

\noindent
{\bf Proof.}
Recall that by (\ref{26.1 B}), we have 
$\vol_g(M)\leq V_0.$ Also, by (\ref{estimate for number in ball 0}),
the $\mu-$volume of a metric ball $B_M(x,\delta_1/6)\subset M$
is bounded from below by $c_3(\delta_1/6)^n$. 

Let
$\{z_1,\dots,z_m\}$ be a  maximal $(\delta_1/3)$-separated subset in $M$.
Then
\ba
m\leq  m_0=1/(c_3(\delta_1/6)^n)=6^n c_3^{-1}\delta_1^{-n}=C_5\delta_1^{-n}.
\ea

Let  $V_k=\{y\in M:\ \dist(y,z_k)< \dist(y,z_j)\hbox{ for all }j\not=k\}$ be the open
Voronoi sets  corresponding to points $z_k$ and
let $W_k$,  $k=1,2,\dots,m$, be such disjoint sets that $V_k\subset W_k\subset \overline{V}_k$ and that
the union of the sets $W_k$  is $M$.
Note that then the balls $B_M(z_j,\delta_1/6)$  are disjoint and
thus
there is $C_6=C_6(n,\Lambda,D,i_0,\rho_{max},\rho_{min})$  so that 
 $$
 \mu(W_k)\geq {\htext c_3} (\delta_1/6)^{-n} \geq 1/(C_6m)
 $$  
 for all  $k=1,2,\dots,m$.
Observe that $\diam(W_k)<\delta_1$, and that if for all $k=1,2,\dots,m$ there is $X_j$, $j\leq N_0$  such that $X_j\in W_k$, then
the set $\{X_j:\ j=1,2,\dots,N_0\}$  is a $\delta_1$-net in $M$.

We can use the classical collectors problem to estimate the probability of the event $A_{m,N_0}$ that
all sets $W_k$ contain at least one point $X_j,$ $j=1,2,\dots,N_0$. 
The tail estimates are used to give a solution for this problem, 
and for the convenience of the reader we give the details of this below
%
(see also \cite{Erdos,Flatto,Newman} for related results).

Let us choose an infinite sequence of i.i.d.\ random variables $X_1,X_2,\dots$
having distribution $\mu$ on $M$
  and let $T$ be the smallest number such that all sets $W_k$, $k=1,2,\dots,m$
contain at least one point $X_j$, $j=1,2,\dots,T$.
 Let $ {E}_i^r$ denote the event that the $i$-th set $W_i$ does not contain any of the first $r$ points $X_1,\dots,X_r$.
Let $m_1=C_6m_0\ge C_6 m $ and $b>1$. Then 
\ba
    \Prob\left [ {E}_i^r \right ] \le \left(1-\frac{1}{C_6m}\right)^r \le \left(1-\frac{1}{m_1}\right)^r \le e^{-r/ {m_1}}.
    \ea
For  $r = \lfloor b m_1 \log m_1\rfloor+1$, we have
$
\Prob\left [ {E}_i^r \right ] \le e^{(-b m_1 \log m_1 ) / m_1} = m_1^{- b}.
$
Then,
\ba
   \Prob\left [ T > b m_1 \log m_1 \right ] = \Prob \left [ \bigcup_{i=1}^m {E}_i^r \right ] \le 
   \frac {m_1}{C_6} \cdot \Prob [ {E}_1^r ] \le \frac 1 {C_6} m_1^{-b + 1}.
   \ea
 Hence, when
 $b= 1+ \frac {\log(2(C_6\theta)^{-1})}{\log m_1}$,
we have
  \ba
  \Prob\left [ T > b m_1 \log m_1 \right ] \le \frac 1 {C_6} m_1^{-b+1}\le \frac \theta 2.
 \ea
 Observe that $m_1= C_6C_5\delta_1^{-n}$ and
  \ba
  b m_1 \log m_1
  = (1+ \frac {\log(2(C_6\theta)^{-1})}{\log (C_6C_5\delta_1^{-n})}) C_6C_5\delta_1^{-n}\log(C_6C_5\delta_1^{-n})
    = (\log(C_6C_5\delta_1^{-n})+ {\log(2(C_6\theta)^{-1})}) C_6C_5\delta_1^{-n}.
  \ea
  Therefore,  when  (\ref{e: N0}) is valid with a suitable $C_3$, we have 
$N_0>  b m_1 \log m_1$. This
 implies
 $
   \Prob(A_{m,N_0})= \Prob(T\leq N_0)\geq  1- \frac \theta 2 .
 $
\hfill \proofbox

\section{The modified $L^2$-norm of the differences of the distance functions}

Let {$y,z\in M$  be (deterministic) points on $M$. 
Denote 
\ba
& &\Phi_y(x)=\Phi(x,y),\quad
\Phi^1_y(x)=\Phi^1(d_M(x,y)).
\ea

\begin{definition}\label{def: ry and A}
For $x\in M$, let $r_x\in C(M)$ we define the distance function 
$$
r_x(y)=d_M(x,y),\quad y\in M.
$$
For $y,z\in M$, let
\beq
k_\Phi(y,z)=\|(r_y-r_z)\Phi^{1/2}_y\Phi^{1/2}_z\|_{L^2(M,d\mu)}^2,\quad 
\Afun{y,z}=\int_M \Phi_y(x)\Phi_z(x)d\mu(x).\hspace{-2cm}
\eeq

\end{definition}

The map $R:M\to L^\infty(M)$, given by $R(x)=r_x$, defines  an isometric embedding $R:M\to R(M)\subset L^\infty(M)$. Below we will
consider this map in different functions spaces. The function $\Afun{y,z}$ measures the relative density of the points $x\in M$ for which the both distances $d_M(x,y)$ and
 $d_M(x,z)$ are non-missing in the data that is given to us. In the next lemmas we analyze these functions.
Recall that $r_1=r_0/2$.

\begin{lemma}\label{lem: geometric observations 1}
 If $d_M(y,z)\leq r_0$ then 
$$
\Afun{y,z}=\int_M \Phi_y(x)\Phi_y(x)d\mu(x)\geq c_4=c_3\phi_1^2r_1^n.
$$
\end{lemma}

\noindent
{\bf Proof.} 
Let $x,y\in M$  be such that $\ell=d_M(y,z)\leq r_1$. Also, let $[yz]=\gamma_{y,\xi}([0,\ell])$, $\xi\in S_yM$ be a length minimizing geodesic from $y$ to $z$. 
Let   $q=\gamma_{y,\xi}(\ell/2)$.
Using properties of $\Phi$ given in (\ref{inequality phi 1}) (see also (\ref{Phi cond 1})-(\ref{Phi cond 2})), we see that for all $x\in B_M(q,r_1)$
we have $x\in B_M(y,r_0)$ and  $x\in B_M(z,r_0)$, and so
$
\Phi_y(x)\Phi_z(x) \geq {}{} \phi_1^2
$.
Observe that 
\beq\label{3H estimate}
\mu(B_M(q,r_1))
\geq
c_3r_1^n.
\eeq
Hence,
\ba
\Afun{y,z}&=&\int_M \Phi_y(x)\Phi_y(x)d\mu(x)\geq
\int_{B_M(q,r_1)} \Phi_y(x)\Phi_y(x)d\mu(x)\geq {}{}\phi_1^2c_3r_1^n.
\ea
\vspace{-5mm}

\hfill \proofbox


Let $X$   have distribution $\mu$.
Let $Y_{X,y}$ be a random variable, taking values in $\{0,1\}$,  that is 1 with probability $\Phi(X,y)$ and
$Y_{X,z}$ be a random variable, taking values in $\{0,1\}$,  that is 1 with probability $\Phi(X,z)$.
Also, let
$\eta,\eta'$ have the zero mean and variance $\sigma^2$,
be such that all $X$,  $\eta$, $\eta'$ are independent random variables.
We assume that under the condition that $X$  is given, random variables $Y_{X,y}$, $Y_{X,z}$,
$\eta$,  and $\eta'$ are independent.

\begin{lemma}\label{lem: samples K}
Let $y,z \in M. $ We have
\beq\label{E one term}
& &\Expec ((|(d_M(y,X)+\eta)-(d_M(z,X)+\eta') |^2-2\sigma^2)
Y_{X,y}Y_{X,z})=
k_\Phi(y,z)\eeq
\end{lemma}

\noindent {{\bf Proof.}
We denote
$ R_y(X)=d_M(y,X)+\eta$ and 
$R_z(X)=d_M(z,X)+\eta'.
$
{Then
\ba
 \Expec  (| R_z(X)- R_y(X)|^2Y_{X,y}Y_{X,z})
&=& \Expec_{\eta,\eta'}\int_M |(d_M(y,x)+\eta)-(d_M(z,x)+{\eta'})|^2\Phi_y(x)\Phi_z(x) d\mu(x)\\
&=&\|(r_y-r_z)\Phi^{1/2}_y\Phi^{1/2}_z\|_{L^2(M,d\mu)}^2+2\sigma^2\Afun{y,z}.
\ea
\vspace{-9mm}

\hfill \proofbox

}}


{\trtext 
\subsection{Deterministic estimates for the rough distance function}

In this subsection, we consider the rough distance function $k_\Phi(y,z)$.

In the study of metric spaces,  Kuratowski observed that the map $R(x)=r_x$
defines an isometric embedding $R:M\to R(M)\subset C(M)$ of the manifold $M$ into
the vector space $C(M)$.
When there are no missing data, that is, $\Phi=1$, the following proposition show that 
the map $\overline R:M\to L^2(M)$, given by $\overline R(x)=r_x$ defines a bi-Lipschitz
embedding $\overline R:M\to \overline R(M)\subset L^2(M)$. Note that here $L^2(M)$
is the space $L^2(M,d\mu)$, where $\mu$ is a probability measure on $M$.


\begin{proposition}\label{prop: L2 norm of distance functions pre}
There is a constant $c_5\in (0,1)$ 
such that
\beq\label{aa inequality 1}
c_5 d_M(y,z)
&\leq &
\|r_y-r_z\|_{L^2(M,d\mu)}\leq d_M(y,z).
\eeq

\end{proposition}

Due to this, we call the map $\overline R:M\to L^2(M)$ the $L^2$-Kuratowski embedding.
The proof of the Proposition \ref{prop: L2 norm of distance functions pre} is 
the special case of the Proposition \ref{prop: L2 norm of distance functions} when $\Phi=1$, claims (i)-(ii), given below.

\begin{proposition}\label{prop: L2 norm of distance functions}

{ 
(i) We have for all $y,z\in M$ the inequality
\beq\label{aa inequality}
\|(r_y-r_z)\Phi^{1/2}_y\Phi^{1/2}_z\|_{L^2(M,d\mu)}\leq \Afun{y,z} d_M(y,z)\leq d_M(y,z).
\eeq

(ii) Let 
$\hat c_4= \frac 14 \min( c_2 H r_1,c_4 ).$
There is $c_5
\leq 1$  such that if $\Afun{y,z}\geq \hat c_4$ then 
\beq\label{ry-rz inequality}
\|(r_y-r_z)\Phi^{1/2}_y\Phi^{1/2}_z\|_{L^2(M,d\mu)}=k_\Phi(y,z)^{1/2}\geq c_5d_M(y,z).
\hspace{-2cm}
\eeq 

(iii) For all $y,z\in M$ satisfying $d_M(y,z)\leq r_1$, 
where $r_1$ is defined in (\ref{r0 and r1 and phi1}), 
we have $
\Afun{y,z}\geq  c_4\ge \hat c_4,$ and so the inequality (\ref{ry-rz inequality}) is valid.

}

\end{proposition}

By this proposition, if $\Afun{y,z}\geq \hat c_4$ then $k_\Phi(y,z)^{1/2}$ approximates $d_M(y,z)$.

\medskip

\noindent
{\bf Proof.}
(i) We have by triangular inequality
\ba
\|(r_y-r_z)\Phi^{\small \tfrac 12}_y\Phi^{\small \tfrac 12}_z\|_{L^2}^2&=&\int_M |d_M(y,\,\cdotp)-d_M(z,\,\cdotp)|^2\Phi_y\Phi_zd\mu
\leq |d_M(y,z)|^2 \Afun{y,z}.
\ea
As $\Afun{y,z}\leq 1$, this proves the inequality (\ref{aa inequality}).

{(ii) To prove  the inequality in (\ref{ry-rz inequality}),
we use the following  (well known) corollary of Toponogov's theorem. Similar kind of formulas are used in 
 Section 4.5 of \cite{Burago}.
However, 
we present the results in the form needed later and give the proof for the convenience of the reader.


\begin{lemma}\label{lem: Topog estimates general}
Let $M$ be a Riemannian manifold with sectional curvature bounded below by $-\Lambda^2$.
Let $x,y,z\in M$ and $\beta=\angle xyz$ be the angle of the length minimizing
curves $[xy]$ and $[yz]$ at $y$.  Assume that $d_M(y,z)\leq \frac 12 d_M(x,y)$
and $d_M(x,y)\leq  \frac 23\min(i_0,\pi/(2\Lambda))$.
Then 
\beq
\label{Si far points new}
 \bigg| d_M(x,z)-(d_M(x,y)-d_M(y,z)\cos\beta)\bigg| \leq  \frac {d_M(y,z)^2}{d_M(x,y)}.
\eeq

\end{lemma}

\noindent
{\bf Proof.} 
Let $\gamma_{y,\xi}([0,\ell])$  be a distance minimizing geodesic from $y$  to $z$, where $|\xi|=1$ and $\ell=d_M(y,z)$. Consider functions
$$
 F(p)=d_M(x,p),\quad f(s)=F(\gamma_{y,\xi}(s)).
$$
Let $\ell_0=\min(i_0,\pi/(2\Lambda))$. Observe that then $d_M(\gamma_{y,\xi}(s),x)\leq \ell_0$ for all $s\in [0,\ell]$.

 The gradient of $F(p)$ at $p\in B(0,\ell_0)$  is equal to the normal vector $\nu$ of the sphere $\Sigma=\p B(x,r)$, where $r=d_M(p,x)$, at the point $p$ and the Hessian of $F$ at $p$ and the shape operator $S(p)$ of the sphere $\Sigma$ have the relation
 Hess$(F)(\xi,\eta)=g(S(p)\xi,\eta)$, $\xi,\eta \in T_pM$, where $g: T_pM\times T_pM\to \R$
 is the quadratic form determined by the metric tensor $g$, see  \cite{Pe}. 
By the standard comparision estimates \cite[Ch. 6, Thm. 2.1]{Pe}, in the space $T_p\Sigma$ we have
\beq
\frac {\Lambda \cosh(\Lambda r)}{\sinh(\Lambda r)} \leq S(p)\leq \frac {\Lambda \cos(\Lambda r)}{\sin(\Lambda r)}.
\eeq
As $\frac d{dt}(\tan(t))=1/\cos^2(t)$, the mean value theorem implies that for $0<s<\pi/(2\Lambda)$ we have
$
\tan(\Lambda s)\geq  \Lambda s,
$
so that $\|S(p)\|\leq 1/ F(p).$ 

Since 
$
\p_s f(s)=g( \nabla F(\gamma_{y,\xi}(s)),\dot \gamma_{y,\xi}(s))
=g(\nu(\gamma_{y,\xi}(s)),\dot \gamma_{y,\xi}(s)),
$
where $\nu(x)=\nabla F(x)$  is the normal of the sphere $\p B(y,s)$ at
the point $x=\gamma_{y,\xi}(s)$. Moreover, $f(0)=d_M(x,y)$ and 
$\p_s f(0)=g(-\dot \gamma_{y,\xi}(s),\dot \gamma_{y,\xi}(s))=-\cos \beta.
$ 
Also,
since 
\ba
\p^2_s f(s)= (\hbox{Hess}\, F)(\dot \gamma(s),\dot \gamma(s))+
g( \nabla F(\gamma(s)), \nabla_{\gamma(s)}\dot \gamma(s))
=g(S(\gamma(s))\dot\gamma(s),\dot \gamma(s)),
\ea
where $\gamma=\gamma_{y,\xi}$,
we have
\ba
|\p^2_s f(s)|\leq \frac {1}{f(s)}\leq  \frac {1}{d_M(x,y)-d_M(y,z)}.
\ea
Hence, using Taylor's series we see that 
$$
\bigg|f(s)-(d_M(x,y)-s\cos \beta)\bigg |\leq \frac{s^2}{2 (d_M(x,y)-d_M(y,z))}\leq \frac{s^2}{d_M(x,y)}.
$$
This proves the claim.
\proofbox

Next we continue the proof of inequality (\ref{ry-rz inequality})
We consider the claim in two cases:

 {\bf Case 1.}    Assume that $d_M(y,z)\ge r_1/16$.
We show that there $c_5'>0$  such that  
\beq\label{aa inequality Phi}
\|r_y-r_z\|_{L^2(M,d\mu)}\geq \|(r_y-r_z)\Phi^{1/2}_y\Phi^{1/2}_z\|_{L^2(M,d\mu)}\geq c_5'd_M(y,z).\hspace{-2cm}
\eeq
 {\bf Proof in Case 1.}  Assume that  $0<d_M(y,z)<r_1/16$. Then, if  $x\in M$  is such that $r_1/2<d_M(x,z)<r_1$, 
we have $d_M(x,y)\geq r_1/4$ and by (\ref{Si far points new}),
\beq
d_M(x,z) 
\nonumber &\leq& d_M(x,y)-d_M(y,z)\cos\beta + \frac 14 d_M(y,z),
\eeq
where $\beta$  is the angle $\angle xyz$.
When $\beta<\pi/4$, this yields 
\ba
d_M(x,y)-d_M(x,z) &\geq &d_M(y,z)\cos\beta - \frac 14 d_M(y,z)
\geq \frac 14 d_M(y,z).
\ea
Thus, let 
\ba
W=\{x\in M;\ r_1/2<d_M(x,z)<r_1,\ \angle xyz <\pi /4\}.
\ea
Then
\ba
\|(r_y-r_z)\Phi^{1/2}_y\Phi^{1/2}_z\|_{L^2(M,d\mu)}^2
&\geq &{}{}\frac 1{16}d (y,z)^2\phi_1^2 \,\vol_\mu(W),
\ea
where by \cite[Cor. 2.4 in Chapter 6.2]{Pe}, see  also (\ref{estimate for number in ball 0}), \eqref{r0 and r1 and phi1},  {\atext there is $c_3'=c_3'(n,\Lambda)>0$ such that}
\ba
\vol_\mu(W)\ge \rho_{min} \vol_M(W)
 \ge
\rho_{min} 4^{-n} \omega_n
{\atext \bigg(\frac {\sin(\Lambda r_1)}{\Lambda r_1 }\bigg)^{n-1}}
(r_1^n-(r_1/2)^n) \geq  
2^{-4n} \rho_{min}{ \atext c_3'} 
r_1^n.
\ea
%
%
%
%
Thus there exists $c_5'$ such that (\ref{aa inequality Phi}) is valid.

{\bf Case 2.} 
Assume that $d_M(y,z)\ge r_1/16$.
Then we show that there $c_5''$  such that  
\beq\label{ry-rz inequality 2}
 \|(r_y-r_z)\Phi^{1/2}_y\Phi^{1/2}_z\|_{L^2(M,d\mu)}=k_\Phi(y,z)^{1/2}\geq c_5''d_M(y,z).
\hspace{-2cm}
\eeq

{\bf Proof in Case 2.} 
The assumption $d_M(y,z)\ge r_1/16 $ and definition of $\hat c_4$ imply that 
$ 4 \hat c_4 c_2^{-1}H ^{-1}\le     r_1$ so that   $d_M(y,z)\ge \frac 14 \hat c_4 c_2^{-1}H^{-1}$. Denote
$d_M(y,z)=\ell$. Then $\ell\geq 2a$, where $a=\frac 18 \hat c_4 c_2^{-1}H^{-1}.$

Assume that $\Afun{y,z}\geq \hat c_4$. Since $d_M(y,x)+d_M(x,z)\geq d_M(y,z)=\ell$,
we see that for all $x\in M$ we have either  $d_M(x,y)\geq \ell/2$ or  $d_M(x,z)\geq \ell/2$.
Let us assume that the latter is true. As $s\mapsto \Phi^1(s)$  is non-increasing and $\Phi$  takes values in $[0,1]$, we have
$
\Phi_y(x)\Phi_z(x)\leq \Phi_z(x) \leq  c_2\Phi^1(d_M(x,z))\leq c_2\Phi^1(\ell/2).
$
This yields that 
\ba
\hat c_4\leq \Afun{y,z}=\int_M \Phi_y(x)\Phi_z(x)d\mu(x)\leq c_2\Phi^1(\ell/2),
\ea
so that $\Phi^1(\ell/2)\geq \hat c_4/c_2 $.

Let $[yz]=\gamma_{y,\xi}([0,\ell])$, $\xi\in S_yM$ be a length minimizing geodesic from $y$ to $z$. 
Let $q=\gamma_{y,\xi}(\ell/2)$,  $p=\gamma_{y,\xi}(\ell/2-a)$,   and $r=a/2$.
When $d_M(x,p)<r$, we have
\ba
d_M(y,x)\leq d_M(y,p)+d_M(p,x)\leq
\ell-\frac a2,\ \ 
d_M(z,x)\geq d_M(z,p)-d_M(p,x)\geq 
\ell+\frac a2,
\ea
so that
$
d_M(z,x)-d_M(y,x)\geq 
a.
$
Recall that $\Phi(x,y)\geq c_1\Phi^1(d_M(x,y))$, $\Phi^1:[0,\infty)\to [0,1]$   and $\|\Phi^1\|_{C^{1}(\R)}\leq H$.
As $\Phi^1_y(q)=\Phi^1_z(q)=\Phi^1(\ell/2)\ge \hat c_4/c_2$ and $B_M(p,r)\subset B_M(q,\frac 32a)$,
we have for all $x\in B_M(p,r)$
\ba
\Phi_y(x)\Phi_z(x)
 \geq c_1^2 (\Phi^1_y(p)-H(a+r))\,(\Phi^1_z(p)-H(a+r))\geq c_1^2 ( c_2^{-1} \hat c_4-\frac 3 2Ha)^2 
\geq  \frac 14 c_1^2 c_2^{-2} \hat c_4^{\,2}.
\ea
Then, as $\ell/D\leq 1$,
\ba
\|(r_y-r_z)\Phi^{1/2}_y\Phi^{1/2}_z\|_{L^2(M,d\mu)}^2&\geq &\int_{B_M(p,r)} |d_M(y,x)-d_M(z,x)|^2 
\Phi_y(x)\Phi_z(x)\,d\mu(x)
\\
&\geq&(c_5'') ^2\ell^2=(c_5'') ^2d_M(y,z)^2,
\ea
where $c_5''= ( 16^{-n-3}   {\rho_{min}}c_1^2   c_2^{-n+4}
\hat c_3 H^{-(n+2)}(\hat c_4)^{n+4} D^{-2})^{1/2}$.
This yields (ii) with 
$c_5=\min (c_5',c_5'')$.

(iii)  As $c_1\le 1$, Lemma \ref{lem: geometric observations 1}
implies that if $d_M(y,z)\leq r_1$, then 
$
\Afun{y,z}\geq c_4\ge \hat c_4.
$
} \hfill \proofbox


\subsection{Probabilistic estimates 
for the rough distance function}

In this subsection, we consider the rough distance function
$k_\Phi(X_j,X_k)^{1/2}$.

Next we determine approximately $k_\Phi(X_j,X_k)=\|(r_{X_j}-r_{X_k})\Phi^{1/2}_{X_j}\Phi^{1/2}_{X_k}\|_{L^2(M,d\mu)}^2$ for $(j,k)\in \I^{(0)}\times \I^{(1)}$ using averaging over the data on
the ``densest net'' $S_2=\{X_j:\ j\in \I^{(2)}\}$.

{\hhtext
For $(j,k)\in I^{(0)}\times I^{(1)}$, we consider the random variables
\beq\label{KN-def B}
K_{jk}&=&\sum_{\hhtext \ell \in I^{(2)}} 
\frac 1{\hhtext N_2} ( |\overline D_{j,\ell }-\overline D_{k,\ell}|^2-{\hhtext 2\sigma^2})Y_{j\ell}Y_{k\ell},\\
%
%
%
%
%
\label{Kjk2 B}
K^L_{jk}&=&\sum_{\hhtext \ell \in I^{(2)}} 
\frac 1{\hhtext N_2} ( \min(|\overline D_{j,\ell }-\overline D_{k,\ell}|^2,L)-{\hhtext 2\sigma^2})Y_{j\ell}Y_{k\ell}.
\eeq
}
We also consider the random variables
\beq\label{Kjk2}
& &A_{jk}=\Afun{X_j,X_k}=\int_M \Phi_{X_j}(x)\Phi_{X_k}(x)d\mu(x),
\quad 
T_{jk}
=\sum_{ \ell \in I^{(2)}}{Y_{j\ell}Y_{k\ell}}  .
\eeq
Roughly speaking, below $T_{jk}$ measures how well $K_{jk}^L$ approximates $k_\Phi(X_j,X_k)$ 
that further approximates
$d_M(X_j,X_k)^2$.

\subsubsection{Probabilistic notations}
To introduce some notations, let us assume for simplicity that the complete probability space $(\Omega,\Sigma,\Prob)$   can be represented as a product
$\Omega=\Omega_1\times \Omega_2\times \Omega_3$  such that $\Prob=\Prob_1\times\Prob_2\times\Prob_3$
and $\omega=(\omega_1,\omega_2,\omega_3)\in \Omega_1\times \Omega_2\times \Omega_3$.
We assume that $\eta=\eta(\omega_1)$, $\eta=\eta'(\omega_2)$ are random variables with variance $\sigma$,
and
 $X=X(\omega_3)$. We assume that $X$ is a random variables having distribution $\mu$. 
Assume that $X$, $\eta$ and $\eta'$  are independent.
 
  For an integrable function $F(\omega_1,\omega_2,\omega_3)=f(\eta(\omega_1),\eta'(\omega_2),X(\omega_3))$
we denote
\beq\label{Expec eta etaprime}
 \Expec_{\eta,\eta'} F= \Expec_{\eta,\eta'}f(\eta,\eta',X) =\int_{\Omega_1}\int_{\Omega_2} F(\omega_1,\omega_2,\omega_3)d\Prob_1(\omega_1)
 d\Prob_2(\omega_2).
 \eeq
 As $X$ is a random variable, we have that  also $\Expec_{\eta,\eta'} f(\eta,\eta',X)$ is a random variable.
  The expectation $ \Expec_{\eta,\eta'}f(\eta,\eta',X)$ over variables ${\eta,\eta'}$ is function of $X$, and thus
 it can be considered as the expectation of $f(\eta,\eta',X)$ under the condition that $X$ is known.
 
Below, we consider   conditional expectations using $\sigma$-algebras.
 Let $\B_X\subset \Sigma$ be
$\sigma$-algebra generated by the random variable $X:\Omega\to \R$,
that is, the $\sigma$-algebra generated by sets $X^{-1}(S)\subset \Omega$, where $S\subset \R$ is an open set,
see \cite[Ch. 5]{Kallenberg}. 
 We recall that $\Expec (F|\mathcal B_X )(\omega)$  is the $\mathcal B_X$-measurable random variable that satisfies
 $$
\int_S \Expec (F|{\mathcal B_X} )(\omega)\,d\Prob(\omega)=\int_S F(\omega)\,d\Prob(\omega)
 $$
 for all sets $S\in \mathcal B_X$. In formula  (\ref{Expec eta etaprime}), $\Expec_{\eta,\eta'} F$  is in  fact equal to the conditional expectation
 $
 \Expec(F|{\mathcal B_X} )=\Expec (F|{\mathcal B_X} )(\omega)$ of the random variable $F$  with respect to the $\sigma$-algebra 
 $\mathcal B_X\subset \Sigma$, that is, we have $(\Expec_{\eta,\eta'} F)(\omega_3)=\Expec (F|{\mathcal B_X} )(\omega)$ with  $\omega=(\omega_1,\omega_2,\omega_3).$
%
As $X$ is a random variable, $\Expec(Z\,| \B_X)$ is  a random variable, too.
Recall also the notation $\Prob(A\,| \B_X)=\Expec(1_A\,| \B_X)$ for
an event $A\in \Sigma$, where $1_A(\omega)$ is the indicator function of the set $A\subset \Omega$.
Below we use several times the fact that
\beq\label{double expectation}
\Expec(\Expec(Z\,| \B_X))=\Expec(Z),\quad \Expec(\Prob(A\,| \B_X))=\Prob(A).
\eeq

Below, we will consider  the $\sigma$-algebra $\B_{j}\subset \Sigma$ generated by the random variable $X_{j}:\Omega\to \R$.
We also consider the
$\sigma$-algebra $\B_{jk}$ generated by the random variables $X_{j}$ and $X_{k}$.

By Lemma \ref{lem: samples K}, the  conditional expectation of $K_{jk}$, under the condition that $X_j$  and $X_k$  are known, satisfies
$
\Expec\bigg(K_{jk}\,\bigg | \B_{jk}\bigg)=k_\Phi(X_j,X_k)
$
where $
k_\Phi(X_j,X_k)=\|(r_{X_j}-r_{X_k})\Phi^{1/2}_{X_j}\Phi^{1/2}_{X_k}\|_{L^2(M,d\mu)}^2
$
is a random variable. 

\subsubsection{Probabilistic estimates for rough distances $K^L_{jk}$ and reliability values $T_{jk}$}

Below use  the following form of Hoeffding's inequality.

\begin{lemma}[Hoeffding's inequality  \cite{Hoeffding}]
Let $Z_1, \dots, Z_N$ be $N$ i.i.d.\ copies of the random variable $Z$ satisfying $0\le Z\le L$, where $L>0$. Then,
 for $\e>0$, we have 
\ba \Prob\left[\left|\frac{1}{N}\left(\sum_{i=1}^N Z_i\right) - \Expec[Z]\right| \leq \e\right] \geq 1 - 2 \exp(-2N\e^2L^{-2}). 
\ea
\end{lemma}

Below, we will show that  $K^L_{jk}$, defined in (\ref{Kjk2 B}), can be considered to be an approximation of $k_\Phi(X_j,X_k)$
which further approximates $d_M(X_j,X_k)^2$ when $A(X_j,X_k)$ is larger than a suitable threshold value.
{\trtext
Let 
\beq\label{e; e3}
 \e_3<\frac 14c_4,\quad b= \frac 12 c_4.
\eeq

 For  $j\in \I^{(0)}$, we consider
 the events $\mathcal E_j^{(1)}\subset \Omega$ and
$
\mathcal E^{(1)}\subset \Omega$, defined by
\beq\label{event Q}
& &\mathcal E^{(1)}_j=\bigg\{\omega\in \Omega\bigg |\ \forall k\in I^{(1)}\ ( \left|\frac{T_{jk}}{N_2} - A_{jk}\right| \leq {\e_3})\bigg\},\quad 
\mathcal E^{(1)}= \bigcap_{ j\in I^{(0)} } \mathcal E^{(1)}_j.
\eeq
Below, {we use a smooth cut-off functions
$$\psi_\rho(t)=\psi_1(t/\rho^2),\quad \beta_1(t)=1-\psi_1(t).$$
where 
 $\psi_1\in C_0^\infty(\R)$ satisfies $\supp(\psi_1)\subset (-2,2)$ and $\psi_1(t)=1$ for $-1 \leq t\leq 1$
 and $0\leq \psi_1(t)\leq 1$ for all $t\in \R,$ $\psi_1(-t)=\psi_1(t)$, $\|\psi_1\|_{C^{1}(\R)}\leq 2$,
 and the function $\psi_1|_{\R_+}$ is non-increasing.
%
%

Note that if
$\beta_1(T_{jk}(bN_2)^{-1})>0$ then
$
\frac{T_{jk}}{bN_2}\ge 1.
$
Also, if $\mathcal E^{(1)}_j$ happens and $A_{jk}\ge  c_4$ then
\beq\label{A.estimate 2}
\frac{T_{jk}}{N_2}\ge A_{jk}-\e_3 \ge  c_4 -\frac 14 c_4 \ge b.
\eeq
Moreover,  if 
$
\frac 1b\,\frac{T_{jk}}{N_2}\ge 1, 
$
then (\ref{A.estimate 2}) implies 
\beq\label{A.estimate}
A_{jk}\ge b-\e_3 \ge  \frac 14 c_4\ge \hat c_4.
\eeq
For $y,z\in M$, let $Y(z,y)$  be a random variable that is 1 with probability $\Phi(y,z)$ and 
0 with probability $1-\Phi(y,z)$. Assume that for $y,z\in M$ the random variables
$Y(y,z)$ are independent.
Let
\beq\label{Kjk2 CC}
T(y,z)
&=&\sum_{ \ell \in I^{(2)}} {Y(y,X_{\ell})Y(z,X_\ell)} .
\eeq
As $\Expec \frac 1{N_2}T(y,z)= \Afun{y,z}$,  Hoeffding's inequality implies
\beq\label{Hoeffding and A Deterministic}
\Prob\left[\left|\frac{T(y,z)}{N_2} -\Afun{y,z}\right| \leq \e_3\right] \geq 1 - 2 \exp(-2N_2\e_3^2). 
\eeq
As $T_{jk}$ and $T(X_j,X_k)$ have the same distributions
and $A_{jk}=A(X_j,X_k)$ for $j\in \I^{(0)}$ and $k\in \I^{(1)}$, inequality (\ref{Hoeffding and A Deterministic}) implies for the conditional probability, under the condition that $X_j$  and $X_k$ are known, that
%
\beq\label{Hoeffding and Ajk}
\Prob\bigg[\left|\frac{T_{jk}}{N_2} - A_{jk}\right| \leq \e_3\,\bigg| \,\B_{jk}\bigg] \geq 1 - 2 \exp(-2N_2\e_3^2). 
\eeq
Thus
we have by (\ref{double expectation})
$\Prob(\mathcal E^{(1)}_j)\geq 1 - 2N_1 \exp(-2N_2\e_3^2).$
Hence,
\beq\label{prob Ew0}
\Prob(\mathcal E^{(1)})\geq 1 - p^{(1)},\quad p^{(1)}=2N_0N_1 \exp(-2N_2\e_3^2).
\eeq

We recall that  
by Lemma \ref{lem: geometric observations 1} and Proposition \ref{prop: L2 norm of distance functions},
\beq\label{ry-rz inequality3}
& &\hspace{15mm}\bigg(d_M(y,z)\leq r_1\implies 
\Afun{y,z}\geq c_4\ge \hat c_4\bigg),\quad \hbox{and }\\ \nonumber
& &\hspace{-25mm}\bigg(\Afun{y,z}\geq \hat c_4
\implies d_M(y,z)\geq \|(r_y-r_z)\Phi^{1/2}_y\Phi^{1/2}_z\|_{L^2(M,d\mu)}=k_\Phi(y,z)^{\frac 12}\geq c_5d_M(y,z)\bigg).
\hspace{-2cm}
\eeq 

}

\begin{lemma}\label{Lemma aux}
Let $L>2\max(D^2,\sigma)$, $\e_2>0$, and
$
\varepsilon(L):= \beta e^{-( L^{1/2}-D)/2}(D^2+6\beta^2)
$
and consider the events $\mathcal E_j^{(2)}\subset \Omega$, ${j\in I^{(0)}}$, and $\mathcal E^{(2)}\subset \Omega$,
\ba
& &\mathcal E_j^{(2)}=\{
\forall k\in I^{(1)}:\ |K^L_{jk} - k_\Phi(X_j,X_k)| \leq \e_2+\varepsilon(L)\},\quad \mathcal E^{(2)}=\bigcap_{j\in I^{(0)}}\mathcal E_j^{(2)}.
\ea
Then
\beq\label{improved2 B CC3}\hspace{1cm}
& &\Prob(\mathcal E^{(2)})
\geq 1 -p^{(2)},\quad p^{(2)}= 2N_0N_1 \exp(-2N_2\e_2^2L^{-2}). 
\eeq
\end{lemma}

\noindent
{\bf Proof.}
Denote $\eta_{jk\ell}= \eta_{j\ell}-\eta_{k\ell}$
and $D_{jk\ell}=d_M(X_j,X_\ell)-d_M(X_k,X_\ell)$.
\btext Then 
\beq\label{double noise}
\Expec \eta_{jk\ell}=0,\quad \Expec \eta_{jk\ell}^2=2\sigma^2,\quad \Expec e^{ |\eta_{jk\ell}|}\leq
(\Expec e^{|\eta_{j\ell}|})(\Expec e^{|\eta_{k\ell}|})
\leq \beta^2.
\eeq
Let $  {\hhtext r}=L^{1/2}$.  Observe that $ |D_{jk\ell}|\leq D$, so that
if $ |D_{jk\ell}+\eta_{jk\ell}|\geq  {\hhtext r}$ then
$|\eta_{jk\ell}|>  {\hhtext r}-D$.
By (\ref{double noise}), $Y=e^{|\eta_{jk\ell}|}$ satisfies
$$
\Prob(|\eta_{jk\ell}|>  {\hhtext r}-D)\leq 
\frac{\Expec ( e^{ |\eta_{jk\ell}|})}{e^{  {\hhtext r}-D}}=\beta^2 e^{-({  {\hhtext r}-D})}.
$$
Thus, using the fact that $\eta_{jk\ell}$ and $D_{jk\ell}$ are independent, we see using Schwartz inequality that
\ba
 \Expec \bigg  ( |(\min((D_{jk\ell}+\eta_{jk\ell})^2,L)-(D_{jk\ell}+\eta_{jk\ell})^2)Y_{j\ell} Y_{k\ell}|  \ \bigg |\,  \B_{jk} \bigg)
& \leq &
\Expec (\chi_{{}_{|\eta_{jk\ell}|>  {\hhtext r}-D}} (D_{jk\ell}+\eta_{jk\ell})^2|\ \B_{jk})\\
&&\hspace{-3cm}\leq (\Prob(|\eta_{jk\ell}|>  {\hhtext r}-D))^{\frac 12}\,(\Expec ((D_{jk\ell}+\eta_{jk\ell})^4|\ \B_{jk}))^{\frac 12}\\
& &\hspace{-3cm}\leq \beta e^{-(  {\hhtext r}-D)/2}(D^2+6\beta^2)\leq\varepsilon(L).
\ea
By Lemma \ref{lem: samples K},
the above shows that
$$
k_{\Phi}^{L}(X_j,X_k):=\Expec ((\min((D_{jk\ell}+\eta_{jk\ell})^2,L)-2\sigma^2)Y_{j\ell}Y_{k\ell}|\ \B_{jk})
$$ 
satisfies
\beq\label{mL-m}
\bigg|k_{\Phi}^{L}(X_j,X_k)-k_\Phi(X_j,X_k)\bigg|\leq \varepsilon(L).
\eeq

By arguing as in (\ref{Hoeffding and A Deterministic})-(\ref{Hoeffding and Ajk}),
we see that Hoeffding's inequality implies
\beq\label{improved2}\hspace{1cm}
\Prob\left[|K^L_{jk} - k_{\Phi}^{L}(X_j,X_k)| \leq \e_2\Bigg| \B_{jk}\right] \geq 1 - 2 \exp(-2N_2\e_2^2L^{-2}). 
\eeq
Using this, \eqref{double expectation} and (\ref{mL-m}), and summing over $j\in I^{(0)}$, we obtain  (\ref{improved2 B CC3}). 
\hfill \proofbox

\section{Determination of the approximate distances in the coarse net}


Next we assume that
\beq\label{req. rho}
\rho\le r_1=r_0/2\leq \phi_0/(4H),\quad u_0= {\phi_1} {\htext c_3}({\rho/4})^n,\quad u_1=u_0/2,\quad u_2=u_0/4.\hspace{-15mm}
\eeq
%

Next we define ${Q}_{j,j'}$ that will turn out to the approximate 
distances $d_M(X_j,X_{j'})$ for points $X_j$ and $X_{j'}$, where $j,j'\in I^{(0)}$, that are sufficiently close to each other.

\begin{definition} Let $\rho\in (0,1)$ satisfy (\ref{req. rho}). For $j,j'\in I^{(0)}$,
let
\beq\label{Q term 1}
{Q}_{j,j'}&=&\frac{ {V}_{j,j'}}{ {W}_{j,j'}},\\
\label{Q term 2}
 {V}_{j,j'}&=& \frac 1{N_1}{\sum_{k\in I^{(1)}}{\beta_1(T_{jk}(bN_2)^{-1} )}\psi_\rho(K_{jk}^{L})}{Y_{j'k}}\overline D_{k,j'},
 \\
 \label{Q term 3}
 {W}_{j,j'}&=&\frac 1{N_1} {\sum_{k\in I^{(1)}}{\beta_1(T_{jk}(bN_2)^{-1} )}\psi_\rho(K_{jk}^{L})}{Y_{j'k}}.
\eeq
 In the case when $ {W}_{j,j'}$  is zero, we define ${Q}_{j,j'}$ to be $D$.
We define for $j,j'\in I^{(0)}$
 \beq \label{Q term 4}
 d^{app}(X_j,X_{j'})=\left\{\begin{array}{cl}{Q}_{j,j'}, &\hbox{ if ${W}_{j,j'}>{u_2}$},\\
 D, &\hbox{otherwise.}\\ \end{array}\right.
 \eeq
\end{definition}

Roughly speaking, above the function 
$\beta_1(T_{jk}(bN_2)^{-1} )$  measures the reliability of the terms $\psi_\rho(K_{jk}^{L})$ in 
formulas \eqref{Q term 2} and \eqref{Q term 3}, 
and ${W}_{j,j'}$ measures reliability of ${Q}_{j,j'}$ in formula \eqref{Q term 4}. 
The numbers $d^{app}(X_j,X_{j'})$ will be the final approximation for the distances $d_M(X_j,X_{j'})$
for all pairs $(X_j,X_{j'})$ of points that are close to each other.
Observe that ${Q}_{j,j'}$ and $ {W}_{j,j'}$ can be computed from the given data.

 
For {technical purposes, 
we define {deterministic (indexed with (d)) and random (indexed with (r)) functions}
\beq\label{W-}
& &W^{(d),-}(y,z)=\int_M 
\beta_1(A(y,x)b^{-1} )\Phi({z,x})\,\psi_{{\htext \rho/2}}(k_{\Phi}(y,x))d\mu(x),\\
\label{W r-}
& &{W^{(r),-}(y,z)= \frac 1{N_1}\sum_{k\in I^{(1)}}
\beta_1(A(y,X_k)b^{-1} )\Phi({z,X_k})\,\psi_{{\htext \rho/2}}(k_{\Phi}(y,X_k)).
}
\eeq
%
%
The motivation behind defining functions $W^{(d),-}$ and $W^{(r),-}$ is that we can use Hoeffding's inequality to estimate how close $W^{(d),-}(X_j,X_{j'})$ and $W^{(r),-}(X_j,X_{j'})$ are when $X_j$ and $X_{j'}$ are known.
Also, we  show that we have $W_{j,j'}\geq W^{(r),-}(X_j,X_{j'})$ with a large probability.

%
%


\begin{lemma}\label{lemma: near points} 
 If  $d_M(x,z)<r_1$, then
\beq\label{W and A large}
& &W^{{(d),-}}(y,z)\geq u_0.
\eeq
Moreover, when we have $\Phi(x,y)\ge c_1\phi_0$
 for all $x,y\in M$,
 the inequality (\ref{W and A large}) holds
for all $x,y\in M$.
\end{lemma}

\noindent{\bf Proof.}
Let  $$w^{-}(x,y,z)=\beta_1(A(y,x)b^{-1} )\Phi({z,x})\,\psi_{{\htext \rho/2}}(k_{\Phi}(y,x))$$  
and  recall that
by (\ref{req. rho}), $\rho\le r_1$.

When $d_M(x,y)<\rho/4$, 
by Lemma \ref{prop: L2 norm of distance functions} (iii),
we have  that $\Afun{x,y}\geq  c_4$. Also, in the case when $\Phi(x,y)\ge c_1\phi_0\geq \phi_1$
 for all $(x,y)\in M\times M$, we have $\Afun{x,y}\ge c_4$.
Thus in both cases, $\beta_1(A(y,x)b^{-1} )=1$.
Moreover, 
 by \ref{prop: L2 norm of distance functions} (i),
we have $k_{\Phi}(y,x)\leq (d_M(y,x))^2<({\htext \rho/2})^2,$ so that  $\psi_{{\htext \rho/2}} (k_{\Phi}(y,x))=1$.


Also, by (\ref{inequality phi 1}), when $d_M(x,z)<r_1=r_0/2$, we have  $\Phi({z,x})\geq \phi_1$.

Thus,
 when $d_M(x,z)<r_1$ or $\Phi(x,y)\ge c_1\phi_0\geq \phi_1$, and we have $d_M(x,y)<\rho/4$, it holds that
$
w_-(x,y,z)\geq \phi_1.
$
Hence we see that
$$
W^{(d),-}(y,z)=\int_M w_-(x,y,z) d\mu(x) \geq \phi_1 \cdotp \mu(B_M(y,\frac \rho 4))\geq \phi_1 \cdotp {\htext c_3}(\rho/4)^n=u_0.
$$
\hfill \proofbox
}

Let us write for $j,j'\in I^{(0)}$
\beq\label{to do 3}
& &{Q}_{j,j'}={Q}_{j,j'}^{1}+{Q}^{2}_{j,j'},
\hbox{ where }
{Q}^{1}_{j,j'}=
\frac { {V}^{1}_{j,j'}}
 { {W}_{j,j'}},\quad {Q}^{2}_{j,j'}=
\frac { {V}^{2}_{j,j'}}
 { {W}_{j,j'}},\\
 \nonumber
& & {V}^{1}_{j,j'}= \frac 1{N_1}{\sum_{k\in I^{(1)}}{\beta_1(T_{jk}(bN_2)^{-1} )}\psi_\rho(K_{jk}^{L})}{Y_{j'k}}d_M(X_k,X_{j'}),
 \\ \nonumber
& & {V}^{2}_{j,j'}=\frac 1{N_1} {\sum_{k\in I^{(1)}}{\beta_1(T_{jk}(bN_2)^{-1} )}\psi_\rho(K_{jk}^{L})}{Y_{j'k}}\eta_{k,j'},
\eeq
and consider the terms ${Q}^{1}_{j,j'}$ and ${Q}^{2}_{j,j'}$ separately.

First we will show that ${Q}_{j,j'}^{2}$  is small with a large probability when $W_{j,j'}\ge u_2$.
To that end, let $h_0>0$ and ${\mathcal E}^{(3)}_{j,j'}\subset \Omega$, $(j,j')\in I^{(0)}⁄\times I^{(0)}$, and ${\mathcal E}^{(3)}\subset \Omega$ be the events
\beq\label{E3}
& &{\mathcal E}^{(3)}_{j,j'}=\{{({W}_{j,j'}\geq {u_2})\implies }(|{Q}_{j,j'}^{2}|\leq h_0)\},\quad
{\mathcal E}^{(3)}=\bigcap_{(j,j')\in I^{(0)}\times I^{(0)}}{\mathcal E}^{(3)}_{j,j'}.  
\eeq

\begin{lemma}\label{Lemma Q^2} 
For any $h_0\in (0,1)$ we have 
\beq\label{to do 4 B}
\Prob({\mathcal E}^{(3)}
)\leq  1-p^{(3)},\quad p^{(3)}=2N_0^2\exp (-e^{-2\beta}N_1 u_2  h_0^2/4).
\eeq
\end{lemma}

\noindent
{\bf Proof.}
{Let us next recall some basic facts:
Let $a=(a_k)_{k=1}^{N_1}$ satisfy $0\leq a_k\leq 1$ and 
\beq\label{def S}
S_{N_1}=
(\sum_{k=1}^{N_1} a_k)^2/(\sum_{k=1}^{N_1} a_k^2),\quad {{Z}}_{N_1}=\sum_{k=1}^{N_1} a_k,
\quad 
V_{N_1}=\frac 1{{{Z}}_{N_1}}\sum_{k=1}^{N_1} a_k\eta_k, \hspace{-1cm} 
\eeq
where $ \eta_k$ are i.i.d.\ variables,  $\Expec\eta_k=0$  
and $\Expec e^{|\eta_k|}\leq \beta$.
Since $0\leq a_k\leq 1$, we have $a_k^2\leq a_k$, so that $\sum_{k=1}^{N_1} a_k^2 \leq \sum_{k=1}^{N_1} a_k,$ and
\beq\label{S and Z}
S_{N_1}=
(\sum_{k=1}^{N_1} a_k)^2/(\sum_{k=1}^{N_1} a_k^2)\geq 
(\sum_{k=1}^{N_1} a_k)^2/(\sum_{k=1}^{N_1} a_k)\geq \sum_{k=1}^{N_1} a_k={{Z}}_{N_1}.
\eeq
Then using Jensen's inequality for the random variable $R=e^{ \eta_k}$
with concave function $h:[0,\infty)\to \R$, $h(s)= s^t$ with $t\in [0,1]$
 (for which the standard  Jensen's inequality reverses), we 
 obtain $\Expec (R^t)\leq (\Expec R)^t$. This yields
$$
\Expec (\exp(t \eta_k))\leq \beta ^t \leq e^{\beta t}
$$
for $t\in [0,1]$. Hence, as $\Expec \eta_k=0$, the moment generating function of $\eta_k$,
$
M(t)=\Expec (\exp(t \eta_k))
$
 satisfies, by considerations involving Taylor series and the mean value theorem,
$
|M_k(t)-1|\leq c't^2,\quad t\in [0,1],
$
where $c'\leq e^{2\beta}$.
Then, by using the independency of random variables $\eta_k$, we see that the moment generating function of $V_{N_1}$  satisfies for $s\in [0,{{Z}}_{N_1}]$
\ba
\Expec \exp({s} V_{N_1})
=\prod_{k=1}^{N_1} M({s} \frac{a_k}{{{Z}}_{N_1}})
\leq \exp (c'\sum_{k=1}^{N_1} ({s} \frac{a_k}{{{Z}}_{N_1}})^2)\leq \exp (e^{2\beta}\frac {{s}^2}{S_{N_1}} )\le \exp (e^{2\beta}\frac {{s}^2}{Z_{N_1}} ),
\ea
where we use the fact that $\sum_{k=1}^{N_1} (a_k/Z_{N_1})^2=1/S_{N_1}$.
Similarly, by considering $-\eta_k$ instead of $\eta_k$, we see that $\Expec \exp(-sV_{N_1})\leq \exp (e^{2\beta}\frac {{s}^2}{S_{N_1}} )\le \exp (e^{2\beta}\frac {{s}^2}{Z_{N_1}} )$.

Then we see that for 
$0<\lambda\le 1$ and $s\in [0,{{Z}}_{N_1}]$,
\beq
\nonumber
\hspace{-1cm}\Prob(|V_{N_1}|> \lambda)&\leq& \Prob(sV_{N_1}>s \lambda)
+ \Prob(sV_{N_1}<-s \lambda) \label{X exp1}
\leq 2\exp (e^{2\beta}s^2Z_{N_1}^{-1} -s\lambda).\hspace{-2cm}
\eeq
When $s=\lambda Z_{N_1}e^{-2\beta}/2 $,
the above implies 
\beq\label{X exp1 B improved}\hspace{-1cm}\Prob(|V_{N_1}|> \lambda)\leq 
2\exp (e^{2\beta}\cdotp 2^{-2}\lambda^2 Z_{N_1}^2e^{-4\beta} \cdotp Z_{N_1}^{-1} -2^{-1}\lambda^2 Z_{N_1}e^{-2\beta})
=2\exp (-\lambda^2 Z_{N_1}e^{-2\beta}/4)
 .\hspace{-2.6cm}
\eeq

%
%
%

{Next we consider a fixed $j,j'\in I^{(0)}$  and let 
$$
a_k={\beta_1(T_{jk}(bN_2)^{-1} )}\psi_\rho(K^L_{jk}){\hhtext Y_{j'k}},
$$ 
 for $k\in I^{(1)}$, $\lambda=h_0$, and let ${{Z}}_{N_1}=N_1{W}_{j,j'}$  be defined analogously to (\ref{def S}). Then we see that  ${Q}_{j,j'}^2=V_{N_1}$, where
 $V_{N_1}$ is be defined analogously to (\ref{def S}) with $\eta_{k}=\eta_{k,j'}$.
  Note that $\eta_{k,j'}$  are independent of variables $T_{jk}$,  {\hhtext  $K^L_{jk}$, and $Y_{j',k}$. Also,
  $V_{N_1}$ is a function of these variables.
  Let $\mathcal B^*$ be the $\sigma-$algebra generated by all
random variables $T_{jk},$  $K^L_{jk}$, and $Y_{j',k}$,  $j,j'\in I^{(0)}$, $k\in I^{(2)}$. 
 As $Z_{N_1}$ is measurable with respect to the $\sigma-$algebra  $\mathcal B^*$, by applying  (\ref{double expectation}) and  (\ref{X exp1 B improved}), we see  
  \ba
\Prob\bigg(({{W}_{j,j'}\geq {u_2})\implies }(|{Q}_{j,j'}^{2}|\leq h_0)\bigg)
&\geq &1-2\exp (-e^{-2\beta}N_1 u_2  h_0^2/4 ).
 \ea
 By doing this analysis for all $j,j'\in I^{(0)}$ and summing up the results, we obtain the claim.
\hfill \proofbox
\medskip

}

Next we analyse ${Q}_{j,j'}^{1}$. 
We assume that $L$ is so large and $\e_2,\e_3$ are so small that 
\beq\label{def h1}
\e_2+\varepsilon(L)\leq \frac 1{100}\rho^2,\quad \e_3\le \frac b {4} u_1.
\eeq

We denote, see \eqref{W-},
\ba
W^{(d),-}_{j,j'}=W^{(d),-}(X_j,X_{j'}),\quad W^{(r),-}_{j,j'}=W^{(r),-}(X_j,X_{j'}),\quad W_{j,j'}=W(X_j,X_{j'}).
\ea
{Let us consider the events ${\mathcal E}^{(4)}_{j,j'}\subset \Omega$ and ${\mathcal E}^{(4)}\subset \Omega$,
\ba
{\mathcal E}^{(4)}_{j,j'}=\bigg\{\omega\in \Omega:\ (W^{(d),-}_{j,j'}\ge u_0)\implies (W^{(r),-}_{j,j'}\geq \frac 12 u_0 )\bigg\},\quad
{\mathcal E}^{(4)}=\bigcap_{(j,j')\in I^{(0)}\times I^{(0)}}{\mathcal E}^{(4)}_{j,j'}.
\ea

Observe that
\beq\label{Wd Wr}
& &\bigg\{\omega\in \Omega:\ |W^{(d),-}(X_j,X_{j'})-W^{(r),-}(X_j,X_{j'})|\leq  \frac 12 u_0 \bigg\}\subset {\mathcal E}^{(4)}_{j,j'}.
\eeq

We see using Hoeffding's inequality that for $y,z\in M$
that 
\ba
\Prob\bigg( |W^{(d),-}(y,z)-W^{(r),-}(y,z)|\leq  \frac 12 u_0 \bigg)\geq 1 - 2 \exp(-2N_1(u_1/2)^2).
\ea
Thus  
we see using (\ref{Wd Wr})  for all $j,j'\in I^{(0)}$ and  \eqref{double expectation}, and summing the results,
\ba
\Prob({\mathcal E}^{(4)})\geq 1 -p^{(4)},\quad  p^{(4)}= 2N_0^2 \exp(-2N_1(u_1/2)^2).
\ea

Below, to shorten notations, let us denote 
$\mathcal E^{(5)}
=\bigcap_{k\in \{1,2,3,4\}}\mathcal E^{(k)}.$
}

\begin{lemma}\label{Lemma Q-Qm new}
Assume the event ${\mathcal E}^{(5)}$ happens. Then, 
if (\ref{def h1}) is valid,
we have   for all $(j,j')\in I^{(0)}\times I^{(0)}$ the implication
 $
({W^{(d),-}_{j,j'}\ge u_0)\implies
(W_{j,j'}\ge u_2}).
$

%
%
%
%
\end{lemma}

\noindent
{\bf Proof.}
As $\mathcal E^{(5)}$ happens, also the event $\mathcal E^{(2)}$ happens.
If  $0\le s\leq \frac  1{20}\rho^2$, then $\psi_\rho(s-\frac 1{10^2}\rho^2)=1$
and $ \psi_{{\htext \rho/2}}(s)=1$.
On the other hand, 
 if $s>\frac 1 {20}\rho^2$, we have  $4(s-\frac 1{10^2}\rho^{2})
>s$. Thus, as the function $\psi_1|_{\R_+}$ is non-increasing,
\ba
\psi_\rho(s-\frac 1{10^2}\rho^2)\ge  \psi_{{\htext \rho/2}}(4(s-\frac 1{10^2}\rho^2))
\ge   \psi_{{\htext \rho/2}}(s)
\ea 
for all $s\ge 0$.
Hence, we have 
\ba
\psi_\rho(s-(\e_2+\varepsilon(L)))&\ge&
\psi_\rho(s-\frac 1{10^2}\rho^2)
\ge  \psi_{{\htext \rho/2}}(s).
\ea
Thus, as ${\mathcal E}^{(2)}$ happens, 
\beq\label{eq: psi K}
 \psi_\rho(K^L_{jk})&\ge &  \psi_\rho(k_{\Phi}(X_j,X_k)-(\e_2+\varepsilon(L)))
\ge  \psi_{{\htext \rho/2}}(k_{\Phi}(X_j,X_k)).
\eeq


%


When  $  \mathcal E^{(5)}$ happens, also $  \mathcal E^{(4)}$ happens, and we have the implication
$
(W^{(d),-}_{j,j'}\ge u_0)\implies
(W^{(r),-}_{j,j'}\ge u_1).
$

Recall that $\|\beta_1\|_{C^{1}(\R)}\leq 2$.
Thus, when $  \mathcal E^{(1)}$ happens, 
%
\beq
 | \beta_1(T_{jk}(bN_2)^{-1} )- 
 \beta_1(A(X_j,X_k)b^{-1} )| \leq 2{b^{-1}}\e_3.
\eeq
Then
\ba
& & | \beta_1(T_{jk}(bN_2)^{-1})  \psi_{{\htext \rho/2}}(k_{\Phi}(X_j,X_k))- 
 \beta_1(A(X_j,X_k)b^{-1} )  \psi_{{\htext \rho/2}}(k_{\Phi}(X_j,X_k))|
  \\
 &\leq& 2{b^{-1}}\e_3\,\cdotp  \psi_{{\htext \rho/2}}(k_{\Phi}(X_j,X_k))
  \leq 2{b^{-1}}\e_3.
\ea
This and (\ref{eq: psi K}) imply 
\beq\nonumber
 \beta_1(T_{jk}(bN_2)^{-1})   \psi_\rho(K^L_{jk})&\ge&
\beta_1(T_{jk}(bN_2)^{-1})   \psi_{{\htext \rho/2}}(k_{\Phi}(X_j,X_k))\\
 &\ge&
 \beta_1(A(X_j,X_k)b^{-1} )  \psi_{{\htext \rho/2}}(k_{\Phi}(X_j,X_k))-2{b^{-1}}\e_3.
\label{semi-deterministic inequality}
\eeq
Recall that  by the assumption of the claim,
 $
\e_3< \frac b 4 u_1
$.
Computing the average of the inequalities (\ref{semi-deterministic inequality}) times $Y_{j',k}$ over $k\in I^{(1)}$, we obtain
\ba
{W_{j,j'}\ge W_{j,j'}^{(r),-}-2{b^{-1}}\e_3\ge W_{j,j'}^{(r),-}-\frac 12 u_1.}
\ea
Thus
 when ${\mathcal E}^{(5)}$  and thus ${\mathcal E}^{(4)}$ happen,
 we have
 \ba
(W^{(d),-}_{j,j'}\ge u_0)\implies (W^{(r),-}_{j,j'}\ge u_1)\implies
(W_{j,j'}\ge u_1-\frac 12 u_1=u_2).
\ea
\vspace{-9mm}

\hfill \proofbox
\medskip

\begin{lemma}\label{Lemma distances}
When the event  ${\mathcal  E}^{(5)}$  happens and (\ref{def h1}) is valid,  it holds for all $j,j'\in I^{(0)}$
that if $W(X_j,X_{j'}) \geq u_2$, then 
\beq
|{Q}_{j,j'}-d_M(X_j,X_{j'})|\leq  \frac 2{c_5}\rho+h_0.
\eeq

\end{lemma}

\noindent
{\bf Proof.}
Below, in this proof assume that 
the event $\mathcal E^{(5)}$ happens. 
%
%
We will first show using  Lemma \ref{Lemma aux}
that if $W(X_j,X_{j'}) \geq u_2$, then 
\beq\label{Q1 - Qd}
|{Q}_{j,j'}^{1}-d_M(X_j,X_{j'})|\leq \frac 2{c_5}\rho.
\eeq
Consider next the indexes $(j,k)\in I^{(0)}\times I^{(1)}$ for which
 $\psi_\rho(K^L_{jk}  )>0$ and $ {\beta_1(T_{jk}(bN_2)^{-1} )})>0$. Then
 $\psi_\rho(K^L_{jk}  )>0$ implies
$
K^L_{jk}<2\rho^2,
$
and hence, as the event $\mathcal E^{(2)}$ happens,
 \beq\label{k large}
k_{\Phi}(X_j,X_{k})<2\rho^2+\e_2+\varepsilon(L).
\eeq
Moreover, if
$ {\beta_1(T_{jk}(bN_2)^{-1} )})>0$ 
then $T_{jk}N_2^{-1}>b$.
As the event $\mathcal E^{(5)}$ happens,
we have 
\ba
A(X_j,X_{k})=
A_{jk}\ge T_{jk}N_2^{-1}-\e_3  \ge b-\e_3 \ge  \frac 14 c_4\ge \hat c_4,
\ea
see Proposition\ \ref{prop: L2 norm of distance functions} and (\ref{e; e3}). Then 
by  
Proposition \ref{prop: L2 norm of distance functions},
\beq\label{ry-rz inequality3 BB}
k_{\Phi}(X_j,X_{k})^{1/2}\geq c_5d_M(X_j,X_{k}),
\hspace{-2cm}
\eeq 
that implies with (\ref{k large}) that
\beq\label{d small}
d_M(X_j,X_{k})\leq \frac 1{c_5}k_{\Phi}(X_j,X_{k})^{1/2}\leq \frac 1{c_5}\sqrt {2\rho^2+\e_2+\varepsilon(L)}.
\eeq

We denote
$
v_{j,j',k}=\beta_1(T_{jk}(bN_2)^{-1} )\psi_\rho(K_{jk}^{L})Y_{j'k}
$
and observe that if $v_{j,j',k}>0$ then the distance $d_M(X_j,X_{k})$  satisfies (\ref{d small}).

Assume next that $W_{j,j'}\ge u_2$.
Then,
we have ${\sum_{k\in I^{(1)}}
 v_{j,j',k}}=W(X_j,X_{j'})\ge u_2>0$  and 
\ba
Q^{1}(X_j,X_{j'})&=&V^{1}(X_j,X_{j'})/W(X_j,X_{j'})
=\frac{\sum_{k\in I^{(1)}}
  d_M(X_k,X_{j'})\,v_{j,j',k}}
{\sum_{k\in I^{(1)}}
 v_{j,j',k}}.
\ea
This means that we can consider $Q^{1}(X_j,X_{j'})$ as 
 a weighted average of distances $d_M(X_k,X_{j'})$.
 Hence, when the event ${\mathcal E}^{(5)}$ happens, we see  that 
 for all $(j,j')\in I^{(0)}\times I^{(0)}$  such that
 $W_{j,j'}\ge u_2$ 
 we have
\ba
|Q^1_{j,j'}-d_M(X_j,X_{j'})|
&=&
\bigg|\frac{\sum_{k\in I^{(1)}}
 (d_M(X_k,X_{j'})-d_M(X_j,X_{j'}))\,v_{j,j',k}}
{\sum_{k\in I^{(1)}}
 v_{j,j',k}}\bigg| \\
& \leq &\frac{\sum_{k\in I^{(1)}}
  d_M(X_k,X_j)\,v_{j,j',k}}
{\sum_{k\in I^{(1)}}
 v_{j,j',k}}
\leq \frac 1{c_5}\sqrt {2\rho^2+\e_2+\varepsilon(L)},
\ea
where we have used  (\ref{d small}) in the last inequality.
By (\ref{def h1}), here
$
\e_2+\varepsilon(L)<\frac 1{100}\rho^2
$
and the  inequality (\ref{Q1 - Qd}) follows.

As we have assumed that  the event $\mathcal E^{(5)}$ happens, also the event $\mathcal E^{(3)}$ happens,
see
(\ref{E3}). Thus  if ${W}_{j,j'}\ge u_2$, then  ${Q}^{2}_{j,j'}<h_0$.
Combining the above, we see that
\ba
|{Q}_{j,j'}-d_M(X_j,X_{j'})|\leq |{Q}^1_{j,j'}-d_M(X_j,X_{j'})|+ |{Q}_{j,j'}^2|\leq  \frac 2{c_5}\rho+h_0.
\ea
\vspace{-.7cm}

\hfill \proofbox
\medskip

The above proposition means that if $W_{j,j'}\ge u_2$, then the number ${Q}_{j,j'}$ approximates $d_M(X_j,X_{j'})$  with a large probability when all parameters are suitably chosen.



Next we consider the proof of Theorem \ref{proposition Q final introduction}.
We use below 
\beq
\label{e: rho} 
\rho=\frac {2\e_1}{c_5},\ 
h_0=\frac { \e_1}2,\  \e_2=\frac{ \rho^2}{200},\  \e_3= \frac b 4 u_1,\ 
L=
4\log^2\bigg(e^{D/2}
 \frac {200\beta (D^2+6\beta^2)}{\rho^{2}}\bigg).\hspace{-19mm}
\eeq
Note that $\rho\le 1$. 
Then $\varepsilon(L) \leq \rho^2/200$,
implying that $\e_2+\varepsilon(L)\leq \rho^2/100.$
Below we assume that $\e_1<\hat \e_1$, where 
$
\hat \e_1=\min(1,8c_5( {\phi_1 {\htext c_3}})^{-1/n}).
$
Then we have  $\e_3<\frac 14c_4$, see (\ref{e; e3}).


Below, see Lemma \ref{Lemma dense net}, let $N_0$  be
\beq\label{e: N0 new}
 N_0= \lfloor 2C_3  \delta_1^{-n} (\log(\delta_1^{-1})+ {\log(\theta^{-1})})\rfloor.
\eeq
Next we consider the probability of $\mathcal E^{(5)}$. 
We see that $\Prob(\mathcal E^{(5)})\geq 1-p^{(5)}$, where 
$p^{(5)}=p^{(1)}+p^{(2)}+p^{(3)}+p^{(4)}$.
Next we consider there probabilities one by one.

\begin{lemma}\label{Lem P3}
There is $C_7>1$ such that we have $p^{(3)}<\theta/8$, when
\beq
\label{e: P3}&& 
 N_1\ge  C_7{\e_1^{-n-2}}(\log(\frac 1{\delta_1})+\log(\frac 1{\theta})).
\eeq
\end{lemma}

\noindent{\bf Proof.}
Below, we use that $t\leq -\log (1-t)$  for $0<t<1$. 
We see that  
$p^{(3)}<\theta/8$ if
\beq
\label{e: N1}N_1\geq R_1=
{\frac{2^{2n+4}e^{2\beta}\log(\frac{16N^2_0}{\theta})}{\phi_1 {\htext c_3}\rho^n h_0^2}
}
=
\frac{2^{2n+8}e^{2\beta}\log(\frac{16 N^2_0}{\theta})}{\phi_1 {\htext c_3} c_5^{2}\rho^n{\rho^2}}.
\eeq


Next we use that $t\log t\leq t^2$, so that for $t>e$ we have
 $\log(t\log t)\leq 2\log t$. Also, recall that $N_0$  is given in (\ref{e: N0 new}).
 Then we see that  
\ba
R_1=\frac{2^{2n+8}e^{2\beta}\log(\frac{16 N^2_0}{\theta})}{\phi_1 {\htext c_3} c_5^{2}\rho^n{\rho^2}}
\leq C_9{\e_1^{-n-2}}(\log(\frac 1{\delta_1})+\log(\frac 1{\theta}))=P_1
\ea
where $C_9$  is suitable.
Thus (\ref{e: N1}) is valid when $N_1\ge P_1$. This yields that claim.
\hfill \proofbox
\medskip

\begin{lemma}\label{Lem P4}
There is $C_{10}>2C_7$ such that we have $p^{(4)}<\theta/8$ when
\beq
\label{e: P4}&& 
 N_1= \lfloor  C_{10}{\e_1^{-2n}}(\log(\frac 1{\delta_1})+\log(\frac 1{\theta}))\rfloor.
\eeq

\end{lemma}

Note that when $N_1$  is chosen as in (\ref{e: P4}), also the inequality (\ref{e: P3}) is valid.
\medskip

\noindent{\bf Proof.}
Next we estimate
$
 p^{(4)}= 2N_0^2 \exp(-2N_1(u_1/2)^2).
 $
We see that  
$p^{(4)}<\theta/8$ if
\beq
\label{e: N4}N_1\ge R_2=
\frac{\log(\frac{16 N^2_0}{\theta})}{2 (u_1/2)^2}=\frac{2\log(\frac{16 N^2_0}{\theta})}{\phi_1^2 {\htext c_3^2}({\rho/4})^{2n}}.
\eeq
Also, we see that
\beq
 \nonumber
R_2
\leq  \lfloor C_{10}{\e_1^{-2n}}(\log(\frac 1{\delta_1})+\log(\frac 1{\theta})) \rfloor =P_2,
\eeq
where $C_{10}$ is suitably chosen. Thus (\ref{e: N4}) is valid when $N_1\ge P_2$.
\hfill \proofbox
\medskip

\begin{lemma}\label{Lem P2}
There is $C_{12}>0$ such that  $p^{(2)}<\theta/8$ when
\beq
\label{e: P2}&& 
 N_2 \ge C_{12}
\e_1^{-4}(\log (\frac 1 {\e_1}))^4
 \bigg(\log (\frac{1}{\theta})+\log (\frac 1 {\delta_1})+\log (\frac 1 {\e_1})\bigg).
\eeq
\end{lemma}

\noindent{\bf Proof.}
Let $N_0$ and $N_1$ are given in (\ref{e: N0 new}) and   (\ref{e: P4}).
We see that  $p^{(2)}\leq \theta/8$ if
 \beq\label{R3 formula}
& & N_2\ge R_3= \frac 12\e_2^{-2}L^{2}  \log(\frac{16 N_0N_1}{\theta}).
 \eeq
We have
\ba
R_3
 \nonumber
& \le& C_{13} \e_1^{-4}L^{2}  \log\bigg(\frac{16}{\theta}\,\cdotp 2C_3  \delta_1^{-n} (\log(\frac 1{\delta_1})+ \log(\frac 1{\theta}))
\,\cdotp C_{10}{\e_1^{-2n}}(\log(\frac 1{\delta_1})+\log(\frac 1{\theta}))\bigg)\\
\nonumber
    & \le&  C_{12}
\e_1^{-4}(\log (\frac 1 {\e_1}))^4
 \bigg(\log (\frac{1}{\theta})+\log (\frac 1 {\delta_1})+\log (\frac 1 {\e_1})\bigg)=P_3 ,
\ea
where $C_{12}$ and $C_{13}$   are suitable. 
Thus (\ref{R3 formula}) is valid when $N_2\ge P_3$. This yields the claim.
\hfill \proofbox
\medskip

%

\begin{lemma}\label{Lem P1}
There is $C_{14}>0$ such that we have $p^{(1)}<\theta/8$ when
\beq
\label{e: P1}&& 
  N_2 \ge C_{14} \e_1^{-2n}   
    \bigg(\log (\frac{1}{\theta})+\log (\frac 1 {\delta_1})+\log (\frac 1 {\e_1})\bigg).
\eeq

\end{lemma}

\noindent{\bf Proof.}
{Using \eqref{e: rho}, we see that the inequality 
\beq\label{extra condition N2}
p^{(1)}=2N_0N_1 \exp(-2N_2\e_3^2)= 2N_0N_1 \exp(-2N_2\frac {\phi_1^2}{2^{4n+8}} {\htext c_3^2}c_4^2\rho^{2n})<\frac 1{8}\theta \hspace{-1cm}
\eeq
is valid when
$N_2\ge R_4=  {2^{4n+7}}{\phi_1^{-2}{\htext c_3^{-2}}c_4^{-2}\rho^{-2n}}
\,  \log(\frac {16 N_0N_1}\theta ).$
We see that
\ba
R_4
\nonumber
&\leq&  \frac   {2^{4n+7}}{\phi_1^2{\htext c_3^2}c_4^2\rho^{2n}}
 \log\bigg(\frac{16}{\theta}\,\cdotp 2C_3  \delta_1^{-n}(\log(\frac 1{\delta_1})+\log(\frac 1{\theta}))
\,\cdotp C_{10}{\e_1^{-2n}}(\log(\frac 1{\delta_1})+\log(\frac 1{\theta}))\bigg)
\\ \nonumber
&\leq& C_{14} \e_1^{-2n}   
    \bigg(\log (\frac{1}{\theta})+\log (\frac 1 {\delta_1})+\log (\frac 1 {\e_1})\bigg)=P_4,
\ea
where $C_{14}$   is suitable. Thus (\ref{extra condition N2}) is valid when $N_2\ge P_4$. This yields the claim.
\hfill \proofbox

Next we prove Theorems 
\ref{thm 2:manifold} and \ref{proposition Q final introduction}.}
\smallskip

\noindent{\bf Proof} (of Theorem \ref{proposition Q final introduction}).
%
We observe that when $\mathcal E^{(5)}$ happens, by Lemma \ref{lemma: near points} and Lemma \ref{Lemma Q-Qm new}, for all $X_j$ and $X_{j'}$ such that $d_M(X_j,X_{j'})<r_1$ we have $W_{j,j'}\ge u_2$.

Let $N_0$ and $N_1$ be given in (\ref{e: N0 new}) and   (\ref{e: P4}).
%
%
%
%
The  conditions (\ref{e: P1})  and (\ref{e: P2}) are valid when we choose a suitable $C_{15}>1$ and
 \beq\label{e: P1 improved}
N_2\ge \lfloor
C_{15} \e_1^{-2n} 
 \bigg(\log^2 (\frac{1}{\theta})+\log^2(\frac 1 {\delta_1})+\log^8 (\frac 1 {\e_1})\bigg)\rfloor 
\eeq
Then
$
p^{(5)}=p^{(1)}+p^{(2)}+p^{(3)}+p^{(4)}\leq \frac 12\theta.
$
This and Lemma \ref{Lemma dense net} prove that with probability $1-\theta$
we have for all $(j,j')\in I^{(0)}\times  I^{(0)}$ that inequality
(\ref{eq: comparison of approximate distances 1})  holds when
$d_M(X_j,X_{j'})<r_1$, and the inequality (\ref{eq: comparison of approximate distances 2}) holds
when $d_M(X_j,X_{j'})\ge r_1$.

In the case when
$\Phi(x,y)\ge c_1\phi_0,$
 for $(x,y)\in M\times M$,
we see that when the events ${\mathcal  E}^{(5)}$  happens, Lemmas \ref{lemma: near points} and \ref{Lemma Q-Qm new} yield that
we will have with probability $1-\theta$ that ${W}_{j,j'}>{u_2}$ for all $(j,j')$.
This implies that
the inequality (\ref{eq: comparison of approximate distances 1}) holds for all pairs $(X_j,X_{j'})$ with
$j,j'\in \{1,2,\dots,N_0\}$.
\hfill \proofbox

\noindent {\bf Proof (of Thm.\ \ref{thm 2:manifold})} Let  $\e_1={{\delta}}^{3/2}$ and {\atext $\delta_1={\atext \Lambda^{2/3}\delta^{1/2}}/20$.
For $N_0,N_1$ given in (\ref{e: N0 new}) and   (\ref{e: P4}),
%
%
\beq
\label{e: N0 B}
  N_0\leq \lfloor 
 C_{16} \delta^{-n/2}  \bigg(\log(\frac{1}{\theta})+\log (\frac 1 {\delta})\bigg) \rfloor ,
\ \ N_1
    \leq   \lfloor C_{17}
\delta^{-3n}
 \bigg(\log^2 (\frac{1}{\theta})+\log^8 (\frac 1 {\delta})\bigg) \rfloor , \hspace{-20mm}
\eeq
with suitable $C_{16}$ and $C_{17}$. Moreover, $N_2$ satisfies \eqref{e: P1 improved} when we choose  a suitable $C_{18}$ and}
\beq 
\label{e: N2 B}
N_2
  \hspace{-2mm}   & =&\hspace{-2mm}  \lfloor C_{18}
\delta^{-3n}
 \bigg(\log^2 (\frac{1}{\theta})+\log^8 (\frac 1 {\delta})\bigg) \rfloor. \hspace{-5mm}
\eeq

Let  $\hat \delta=\e_1=\delta^{3/2}$ and $\hat r=(\hat \delta/\Lambda^2)^{1/3}={\atext \Lambda^{2/3}\delta^{1/2}}$.
Then by
Theorem \ref{proposition Q final introduction}, with probability $1-\theta$ the set $\mathcal X=\{X_j:\ j=1,2,\dots,N_0\}$
is a $\delta_1$-dense subset of $M$ and the approximate distances $\tilde d(X_j,X_{j'})=
 d^{app}(X_j,X_{j'})$, $j,j'\in \{1,2,\dots  N_0\}$, see (\ref{Q term 4}),
satisfy the conditions given in Proposition \ref{prop:improved density} in Appendix A. Thus with probability $1-\theta$
 can apply
Proposition \ref{prop:improved density} (see also \cite[Corollary 1.10]{FIKLN})
 with  $\hat \delta=\delta^{3/2}$ and $\hat r=(\hat \delta/\Lambda^2)^{1/3}$ to construct a Riemannian manifold $(M^*,g^*)$
 that approximates the original manifold $(M,g)$ so that the  claims (1)-(3) in Theorem \ref{thm 2:manifold} are satisfied.
\hfill \proofbox

\section*{Appendix A: Reconstruction of a manifold with a small deterministic errors }
{\atext Here, we give results on  the reconstruction of a Riemannian manifold
when one is given distances with small deterministic errors. The following result is an
improvement of Corollary 1.10 in \cite{FIKLN}.

\begin{proposition} 
\label{prop:improved density}
There are  $C_n'>0$, depending on $n$, and $c_1'(n,K)>0$, depending on $n,K$, such that the following holds:
Let $0<\hat \delta<c_1'(n,K)$,  $\hat r= (\hat \delta/K)^{1/3}$ and  $M$  be a compact $n$-dimensional manifold 
with $|\Sec(M)|\leq  K$ and $\inj(M)>2\hat r$. Let 
  $\mathcal X=\{x_j\}_{j=1}^N$
be an $\hat r/20$-dense subset of $M$.
 Moreover, let
%
   $\tilde d \co
 \mathcal X\times \mathcal X\to \R_+\cup\{0\}$
be {an approximate local distance function} that satisfies 
\beq
\label{delta-condition 0 B}
& & |\tilde d(x,y)-d_M(x,y)|\leq \hat \delta,\quad \hbox{if } d_M(x,y)< \hat r,\\
\label{delta-condition 1 B}
& & \tilde d(x,y)>\hat r-\hat\de,\quad \hbox{if } d_M(x,y)\ge \hat r.
\eeq
%
%
%
%

Then, given the values $\tilde d(x_j,x_k)$, $j,k=1,2,\dots,N$, one can 
 construct a compact $n$-dimensional Riemannian manifold $(M^*,g^*)$
such that:

\begin{enumerate}
\item 
There is a diffeomorphism $F:M^*\to M$ satisfying 
\ba
\frac 1L\leq \frac{d_{M}(F(x),F(y))}{d_{M^*}(x,y)}\leq L,\quad \hbox{for }x,y\in M^*,\ L=1+C_n'K^{1/3}\hat \delta\,{}^{2/3}.
\ea

\item   $|\Sec(M^*)|\le C_n'K$.

\item The injectivity radius $\inj(M^*)$ of $M^*$  
satisfies
$$
\inj(M^*)\ge \min\{ (C_n'K)^{-1/2}, (1-C_n'K^{1/3}\hat \delta\,{}^{2/3})\inj(M)\} .
$$
\end{enumerate}
\end{proposition}}

\begin{proof}
{\atext A  result similar to the claim is proven in  \cite[Corollary 1.10]{FIKLN} under the assumption that the 
set $\mathcal X$ is a $\hat\delta$-dense subset of $M$, instead of $\hat r/20$-dense as it is assumed in the claim. 
Moreover, by \cite[Corollary 1.10]{FIKLN}, it is enough 
 to construct numbers $\tilde D_{j,k}$, $j,k=1,2,\dots,\tilde N$, such that the following is true:
 %
%
%
%
There a $\hat\delta$-net $\mathcal Y=\{y_j:\ j=1,2,\dots,\tilde N\}\subset M$
such that the conditions  \eqref{delta-condition 0 B} and \eqref{delta-condition 1 B}
are valid for the function $\tilde d'\co \mathcal Y\times\mathcal Y\to \R_+\cup\{0\}$ defined by 
$\tilde d'(y_j,y_k)=\tilde D_{j,k}$.

Next we construct
 the required
   $\hat\delta$-net $\mathcal Y\subset M$
and an approximate distance function $\tilde d'$ on
$\mathcal Y\times\mathcal Y$.}
%
%
%
%
%
%
%
We assume that $c_1'(n,K)$ is chosen so small that
$\hat\delta/\hat r = K\hat r^2 < \frac1{150}$.

For $p\in M$ we denote by $E_p$ the restriction of
the Riemannian exponential map $\exp_p$ 
to the $\hat r$-ball in $T_pM$ centered at the origin.
This restriction is a diffeomorphism onto the $\hat r$-ball centered at $p$ in~$M$.
It distorts distances by at most $\frac12\hat\delta$,
namely for all $u,v\in T_pM$ such that $|u|,|v|<\hat r$ we have
\beq\label{exp distortion}
  \bigl|d_M(E_p(u),E_p(v))-|u-v|\bigr| < \tfrac12 K\hat r^3 =\tfrac12\hat\delta .
\eeq
This inequality holds as long as $\inj(M)>2\hat r$ and $K\hat r^2<\pi/2$,
see \cite[Section 4]{FIKLN} for a proof.

For every $p\in\mathcal X$, define
$X_p = \{ x\in \mathcal X : \tilde d(p,x) < \hat r/6-\hat\delta \}$.
By \eqref{delta-condition 0 B} and \eqref{delta-condition 1 B},
$X_p$ is contained in the $\hat r/6$-neighborhood of $p$.
Define $\tilde X_p=E_p^{-1}(X_p)$,
let $\tilde V_p$ be the convex hull of $\tilde X_p$ in $T_pM$,
and $V_p=E_p(\tilde V_p)$.
Since $\mathcal X$ is $\hat r/20$-dense in $M$,
\eqref{delta-condition 0 B} and \eqref{exp distortion} imply that
for every $u\in T_pM$ such that $|u|<\hat r/6-\hat r/20-2\hat\delta$
there exists $v\in\tilde X_p$ such that $|u-v|<\hat r/20+\hat\delta/2$.
This implies that $\tilde V_p$ contains the ball of radius
$\hat r/6-2\hat r/20-2\hat\delta-\hat\delta/2>\hat r/20$ centered
at the origin. 
(Here we use the assumption that $\hat\delta/\hat r<\frac1{150}$).
Hence $V_p$ contains the $\hat r/20$-ball centered at~$p$
and therefore $\bigcup_{p\in\mathcal X}V_p=M$.

We represent points of $\tilde V_p$ as linear combinations
of points of $\tilde X_p$ as follows.
Let $X_p^n$ be the set of all $n$-tuples of points of $X_p$
and $\Delta^n$ the standard coordinate simplex in $\R^n$:
$$
\Delta^n = \{ (t_1,\dots,t_n)\in\R^n : t_1,\dots,t_n\ge 0, \sum t_i\le 1 \}.
$$
For $\alpha=(a_1,\dots,a_n)\in X_p^n$ and
$\tau=(t_1,\dots,t_n)\in\Delta^n$, let
$$
 S_p(\alpha,\tau) = \sum t_i E_p^{-1}(a_i) .
$$
This defines a map $S_p\co X_p^n\times\Delta^n\to T_pM$.
For a fixed $\alpha=(a_1,\dots,a_n)\in X_p^n$,
the range of $S_p(\alpha,\cdot)$ is a (possibly degenerate)
affine simplex in $T_pM$ with vertices 
$0, E_p^{-1}(a_1),\dots,E_p^{-1}(a_n)$.
Since $0\in\tilde X_p$, the union of all such affine simplices
is precisely the convex hull of $\tilde X_p$.
Thus $S_p(X_p^n\times\Delta^n) = \tilde V_p$.

Fix an $\ep'$-dense finite set $\Sigma\subset\Delta^n$,
where $\ep'=\hat\delta/(3\hat r\sqrt n)$,
and define
$
 Y_p = E_p(S_p(X_p^n\times\Sigma)) \subset M .
$
Since $\tilde V_p$ is contained in the $\hat r/6$-ball,
$S_p(\alpha,\cdot)$ is Lipschitz with Lipschitz constant $\hat r\sqrt n/6$.
Therefore $S_p(X_p^n\times\Sigma)$ is $\hat\delta/2$-dense in $\tilde V_p$.
Hence, by \eqref{exp distortion}, $Y_p$ is $\hat\delta$-dense
in~$V_p$.

Now define  $\mathcal Y\subset M$ by
$\mathcal Y = \bigcup_{p\in \mathcal X} Y_p $.
Since the sets $V_p$ cover $M$
and $Y_p$ is $\hat\delta$-dense in $V_p$ for each~$p$,
$\mathcal Y$ is a $\hat\delta$-net in $M$.
The points of $\mathcal Y$ are indexed by triples $(p,\alpha,\tau)$
where $p\in\mathcal X$, $\alpha\in X_p^n$, $\tau\in\Sigma$.
This index set can be enumerated algorithmically using the known data.

Our first goal is to compute approximate {\em squared} distances
$Q(x,y)$ between sufficiently close pairs of points $x,y\in\mathcal Y$.
Fix $p,q\in\mathcal X$ such that $\tilde d(p,q) < 2\hat r/3-2\hat\delta$
(the case $p=q$ is not excluded).
By \eqref{delta-condition 0 B} and \eqref{delta-condition 1 B} we have
$d_M(p,q)<2\hat r/3-\hat\delta$.
Hence, by the triangle inequality,
$ d_M(x,y) < \hat r-\hat\delta $
for all $x\in V_p$ and $y\in V_q$,
In particular $V_q$ is contained in the range of~$E_p$.
By \eqref{exp distortion},
\beq\label{exp square distortion}
 \left| d_M(x,y)^2- |E_p^{-1}(x)-E_p^{-1}(y)|^2 \right| < \hat\delta\hat r
 \qquad\text{for all $x\in V_p$ and $y\in V_q$}.
\eeq

We compute the values $Q(x,y)$
for all $x\in X_p\cup Y_p$ and $y\in X_q\cup Y_q$
in several steps.
First consider $x\in X_p$ and $y\in X_q$.
In this case we simply define $Q(x,y)=\tilde d(x,y)^2$.
Then, by \eqref{delta-condition 0 B},
\beq\label{Q estimate 2}
| Q(x,y) - d_M(x,y)^2 | < 2\hat\delta\hat r .
\eeq
Hence, by \eqref{exp square distortion},
\beq\label{Q estimate 3}
 \left| Q(x,y) - |E_p^{-1}(x)-E_p^{-1}(y)|^2 \right| < 3\hat\delta\hat r ,
 \qquad x\in X_p, \ y\in X_q .
\eeq

Now consider $x\in Y_p$ and $y\in X_q$.
By the definition of $Y_p$ we have $x=E_p(S_p(\alpha,\tau))$ 
for some $\alpha=(a_1,\dots,a_n)\in X_p^n$
and $\tau=(t_1,\dots,t_n)\in\Sigma$.
We define $Q(x,y)$ using the values of $Q$
that we have from the previous step.
Introduce the following notation: $a_0=p$, $t_0=1-\sum_{i=1}^n t_i$,
$v_i=E_p^{-1}(a_i)$ for $i=0,\dots,n$
(in particular $v_0=0$),
$v=E_p^{-1}(x)=\sum t_i v_i$, and $w=E_p^{-1}(y)$.
In this notation, $ v-w = \sum_{i=0}^n t_i (v_i-w) $,
hence
$$
 |v-w|^2 = \sum_{0\le i,j\le n} t_it_j \langle v_i-w,v_j-w\rangle
 = \frac12 \sum_{0\le i,j\le n} t_it_j (|v_i-w|^2+|v_j-w|^2-|v_i-v_j|^2)
$$
where $\langle\,\cdotp,\,\cdotp\rangle$ is the scalar product in $T_pM$.
With this identity in mind, we define
\beq\label{Q definition}
 Q(x,y) = \frac12\sum_{0\le i,j\le n}
  t_it_j (Q(a_i,y)+Q(a_j,y)-Q(a_i,a_j)).
\eeq
Since $y\in X_q$ and $a_i,a_j\in X_p$, the
values $Q(a_i,y)$, $Q(a_k,y)$ and $Q(a_i,a_j)$
are defined in the previous step
(in the case of $Q(a_i,a_j)$, this is the previous step with $q=p$).
Since $\sum_{i=0}^n t_i=1$,
applying \eqref{Q estimate 3} to the terms
in the right-hand side of \eqref{Q definition}
yields that
\beq\label{Q estimate 5}
 | Q(x,y) - |E_p^{-1}(x)-E_p^{-1}(y)|^2 | = | Q(x,y) - |v-w|^2 |
 < \tfrac 92\hat\delta\hat r < 5\hat\delta\hat r .
\eeq
Therefore, by \eqref{exp square distortion},
\beq\label{Q estimate 6}
 | Q(x,y) - d_M(x,y)^2 | < 6\hat\delta\hat r .
\eeq
We now have the values $Q(x,y)$ satisfying \eqref{Q estimate 6}
for all $x\in Y_p$ and $y\in X_q$.
Exchanging the roles of $p$ and $q$ we similarly find $Q(x,y)$
for all $x\in X_p$ and $y\in Y_q$.

Finally, consider $x\in Y_p$ and $y\in Y_q$. Again, let
$x=E_p(S_p(\alpha,\tau))$ where $\alpha=(a_1,\dots,a_n)\in X_p^n$
and $\tau=(t_1,\dots,t_n)\in\Sigma$, and define $a_0=p$
and $t_0=1-\sum_{i=1}^n t_i$.  From the previous steps
we already have values $Q(a_i,y)$ and $Q(a_i,a_j)$.
Therefore we can define $Q(x,y)$ by the same
formula \eqref{Q definition}.
Then, starting from \eqref{Q estimate 6} instead of \eqref{Q estimate 2},
we obtain the same estimates as above but with different constants:
\eqref{Q estimate 3} with $7\hat\delta\hat r$ in the right-hand side,
\eqref{Q estimate 5} with $11\hat\delta\hat r$ in the right-hand side,
and finally \eqref{Q estimate 6} with $12\hat\delta\hat r$ in the right-hand side:
\beq\label{Q estimate 12}
 | Q(x,y) - d_M(x,y)^2 | < 12\hat\delta\hat r ,
 \qquad x\in Y_p, \ y\in Y_q .
\eeq

Now one might take the square root of $Q(x,y)$
as an approximate distance between $x$ and $y$;
however this approximation is not good enough.
For a better one, we use an algorithm described in
\cite[\S2.4]{FIKLN} to construct a map $F\co Y_p\cup Y_q\to\R^n$
that preserves distances up to an error $O(\hat\delta)$.
Let us outline how the algorithm works in the present set-up.

First define an approximate scalar product $P(x,y)$ for all pairs
$x,y\in Y_p\cup Y_q$ by
$$
 P(x,y) = \tfrac12 (Q(p,x) + Q(p,y) - Q(x,y) ) .
$$
By \eqref{exp square distortion}, \eqref{Q estimate 12}
and the Euclidean identity
$ \langle u,v\rangle = \tfrac12 (|u|^2 + |v|^2 - |u-v|^2 ) $
for $u,v\in T_pM$,
this approximates the scalar product
of $E_p^{-1}(x)$ and $E_p^{-1}(y)$ in $T_pM$:
\beq\label{scalar product estimate}
 |P(x,y) - \langle E_p^{-1}(x),E_p^{-1}(y) \rangle| < 20 \hat\delta\hat r.
\eeq
Then, since $E_p^{-1}(X_p)$ is a $\hat\delta/2$-net
in $E_p^{-1}(V_p)$ and the latter contains the $\hat r/6$-ball
centered at the origin, we can find points $a_1,\dots,a_n\in Y_p$
such that the vectors $v_i:=E_p^{-1}(a_i)$, $i=1,\dots,n$,
approximate an orthonormal basis of $T_pM$
rescaled by the factor $\hat r/6$:
\beq\label{almost orthogonal}
 |  (\hat r/6)^{-2} \langle v_i, v_j \rangle -\delta_{ij} |
 < C_1 \hat\delta/\hat r,
 \qquad 1\le i,j\le n,
\eeq
where $\delta_{ij}$ is the Kronecker delta and
$C_1=C_1(n)>0$ is a suitable constant.
(A straightforward modification of the algorithm
from \cite[\S2.4]{FIKLN} can be used to find such points efficiently).

The inequalities \eqref{almost orthogonal} imply that
the linear map $L\co T_pM\to\R^n$ defined by
$$
L(v)=(\hat r/6)^{-1}(\langle v,v_1\rangle,\dots,\langle v,v_n\rangle)
$$
is $(C_2\hat\delta/\hat r)$-close in the operator norm
to a linear isometry from $T_pM$ to $\R^n$
for some constant $C_2=C_2(n)>1$, see \cite[Lemma 2.6]{FIKLN}.
Hence $L$ distorts distances within the
$\hat r$-ball by at most $2C_2\hat\delta$.
We approximate $L\circ E_p^{-1}$ by a map
$F\co Y_p\cup Y_q\to\R^n$ defined by
$$
 F(x) = (\hat r/6)^{-1}(P(x,a_1),\dots,P(x,a_n)), \qquad x\in Y_p\cup Y_q,
$$
and compute $\tilde d'(x,y)=|F(x)-F(y)|$ for all $x,y\in Y_p\cap Y_q$.
By \eqref{scalar product estimate} we have
$
 |F(x)-L(E_p^{-1}(x))| < 120\sqrt n\,\hat\delta
$
for all $x,y\in Y_p\cap Y_q$.
Hence, by \eqref{exp distortion} and the above mentioned property of $L$,
\beq\label{tilde d estimate}
 |\tilde d'(x,y) - d_M(x,y)| < C_4\hat\delta
\eeq
where $C_4=2C_2+120\sqrt n +1$.

The domain of the function $\tilde d'$ defined by the above procedure
includes all pairs $x,y\in\mathcal Y$ such that $d_M(x,y)<\hat r/4$.
Indeed, if $d_M(x,y)<\hat r/4$ and $p,q\in X$ are such that
$x\in Y_p$, $y\in Y_q$,
then by the triangle inequality we have
$d_M(p,q)<\hat r/4+2\hat r/6<2\hat r/3-3\hat\delta$,
and hence, by \eqref{delta-condition 0 B},
$\tilde d(p,q)<2\hat r/3-2\hat\delta$.
Thus for any such pair $x,y$ the value $\tilde d'(x,y)$
is defined and satisfies \eqref{tilde d estimate}.

To finish the construction, set $\tilde d'(x,y)=\hat r$
for all remaining pairs $x,y\in\mathcal Y$.
Now the function $\tilde d'$ is defined on $\mathcal Y\times\mathcal Y$
and it satisfies the assumptions  \cite[Corollary 1.10]{FIKLN}
for $\hat r'=\hat r/4$ in place of $\hat r$,
$\hat\delta'=C_4\hat\delta$ in place of $\hat\delta$,
and  $K'=\hat\delta'/(\hat r')^3=2^6C_4 K$ in place of~$K$.
Applying  \cite[Corollary 1.10]{FIKLN} 
with these modified parameters finishes the proof
of Proposition \ref{prop:improved density}.
\end{proof}

{\atext The constructions in Prop.\ \ref{prop:improved density} and  \cite[Corollary 1.10]{FIKLN} are algorithmic, for  the details, see \cite{FIKLN}.}

\vfill

\end{document}